\documentclass[11pt]{article}

\title{Stochastic block models with many communities and the Kesten--Stigum bound}

\author{Byron Chin\thanks{Department of Mathematics, Massachusetts Institute of Technology. Email: \textup{\tt byronc@mit.edu}} \and Elchanan Mossel\thanks{Department of Mathematics, Massachusetts Institute of Technology. Email: \textup{\tt elmos@mit.edu}} \and Youngtak Sohn\thanks{Division of Applied Mathematics, Brown University. Email: \textup{\tt youngtak\_sohn@brown.edu}} \and Alexander S. Wein\thanks{Department of Mathematics, University of California, Davis. Email: \textup{\tt aswein@ucdavis.edu}}}

\date{}

\usepackage{arxiv}
\numberwithin{equation}{section}

\begin{document}

\maketitle

\begin{abstract}
  We study the inference of communities in stochastic block models with a growing number of communities. For block models with $n$ vertices and a fixed number of communities $q$, it was predicted in Decelle~et~al.~(2011) that there are computationally efficient algorithms for recovering the communities above the Kesten--Stigum (KS) bound and that efficient recovery is impossible below the KS bound. This conjecture  has since stimulated a lot of interest, with the achievability side proven in a line of research that culminated in the work of Abbe and Sandon~(2018). Conversely,  recent work provides evidence for the hardness part using the low-degree paradigm.

  In this paper we investigate community recovery in the regime $q=q_n \to \infty$ as $n\to\infty$ where no such predictions exist. We show that efficient inference of communities remains possible above the KS bound. Furthermore, we show that recovery of block models is low-degree hard below the KS bound when the number of communities satisfies $q\ll \sqrt{n}$. Perhaps surprisingly, we find that when $q \gg \sqrt{n}$, there is an efficient algorithm based on non-backtracking walks for recovery even below the KS bound. We identify a new threshold and ask if it is the threshold for efficient recovery in this regime. Finally, we show that detection is easy and identify (up to a constant) the information-theoretic threshold for community recovery as the number of communities $q$ diverges. 

  Our low-degree hardness results also naturally have consequences for graphon estimation, improving results of Luo and Gao (2024).\footnote{Accepted for presentation at the Conference on Learning Theory (COLT) 2025}
\end{abstract}

\section{Introduction}

The \emph{stochastic block model}, or simply the \emph{block model}, is a canonical model for networks with community structure, and has been studied extensively in the social sciences, statistics and computer science \cite{HLL:83, DF:89, SN:97, CK:01, M:01, BC:09, C:10, RCY:11}. Work in statistical physics initiated the study of inference of \emph{sparse} block models, where the average degree is constant \cite{DKMZ:11}. In this sparse regime, the symmetric block model is parameterized by the number of vertices $n$, the number of communities $q$, and two parameters $a,b > 0$. The vertices are partitioned into communities of equal size, and edges are formed independently with probability $a/n$ between vertices within the same community and probability $b/n$ between vertices in different communities. A useful re-parametrization is letting $d\equiv \frac{a+(q-1)b}{q}$ denote the average degree and $\lambda\equiv \frac{a-b}{a+(q-1)b}$ measure the strength of the signal of within community versus between community edges.

The main inference task for the block model is to extract some information about the community assignments of the vertices, given a single observed graph drawn from the distribution. We use the vague term ``some information'' here due to $(i)$ inherent symmetries in the model that only allow us to recover communities up to relabeling, and $(ii)$ the sparsity of the model, which, even up to these symmetries, allows for recovery only up to constant correlation with the true communities. The formal definition of the recovery notion we consider, commonly referred to as \emph{weak recovery}, is provided in Definition~\ref{def: Align} below.

The paper of Decelle, Krzakala, Moore, and Zdeborov\'a~\cite{DKMZ:11} conjectured that the signal-to-noise ratio (SNR) $d\lambda^2$ governs the feasibility of computationally efficient (polynomial-time) recovery in block models. Their conjecture, based on deep heuristics from statistical physics, states that weak recovery is computationally feasible when $d \lambda^2 > 1$, but infeasible when $d \lambda^2 < 1$. The critical threshold $d\lambda^2=1$ is known as the \emph{Kesten--Stigum (KS) bound}, which was first pinpointed in the study of multi-type branching processes by Kesten and Stigum~\cite{KS:66}. A long line of work has since studied this threshold in various contexts including the reconstruction on trees problem~\cite{BRZ:95, EKPS:00, MP:03, JM:04} and the phylogenetic reconstruction problem~\cite{M:04, DMR:06, MRS:11, RS:17}. 

For the case of two communities $q=2$, the conjecture of Decelle \emph{et al.}~\cite{DKMZ:11} was rigorously proven in a series of papers, where spectral methods or algorithms based on non-backtracking walks succeed above the KS bound \cite{M:14, MNS:18, BLM:15} while recovery is information-theoretically impossible below it~\cite{MNS:15}. For $q > 2$, the rigorous understanding is less complete. Abbe and Sandon~\cite{AS:18} proved the achievability of weak recovery above the KS bound, and Mossel, Sly, and Sohn~\cite{MSS:22} established the information-theoretic impossibility of weak recovery below the KS bound for $q\leq 4$ and large values of $d$.  However, when $q\geq 5$, information-theoretic recovery is possible below the KS bound~\cite{BMNN:16, AS:15a}. Despite this, the prevailing conjecture remains that the KS bound represents the \emph{computational} threshold for weak recovery for all values of \emph{fixed} $q$~\cite{AS:18}.

While the setting of a growing number of communities has also been studied previously (e.g.,~\cite{CX:16}), there remains a large gap in our understanding of the sparse regime. Notably the powerful methods from statistical physics that give predictions for many such inference questions are believed to be inapplicable in this regime. A primary motivation for our work is to present a picture of the computational phases in this regime where statistical physics methods have not yet made predictions to the best of our knowledge. 

In this paper, we study efficient recovery in the sparse stochastic block model with $q = \omega(1)$ communities. We will often use the parametrization $q = n^\chi$ for some exponent $\chi \in (0,1)$. The central question is:
\begin{quotation}
    \centering{\emph{Does the KS bound remain the threshold for efficient recovery in stochastic block models with a growing number of communities?}}
\end{quotation}
We provide rigorous evidence that the conjectural picture for weak recovery with a finite number of communities persists for $\chi < 1/2$. More specifically, we prove that any $n^{\Omega(1)}$-degree polynomial in the adjacency matrix of the observed graph cannot achieve weak recovery below the KS bound when $\chi<1/2$. Low-degree polynomial estimators have a strong record of tracking computational hardness in average-case problems (for more discussion, see Section~\ref{subsec:low:degree} below); therefore, our result provides strong evidence for the computational intractability of weak recovery below the KS bound when $\chi<1/2$. Moreover, we establish efficient recovery above the KS bound when $q \to \infty$.

Perhaps surprisingly, we establish that for $\chi>1/2$, efficient recovery is achievable using non-backtracking walks even below the usual KS bound $d\la^2=1$. Thus, there is a transition in feasibility of weak recovery with respect to the KS bound at $q \asymp \sqrt{n}$. Our analysis leads us to propose $d\lambda^{1/\chi}$ as the appropriate SNR to determine the threshold for recovery in this regime.

\subsection{Main results}
Formally, the \textit{stochastic block model} is a random graph $G \sim \SBM(n, q, d, \lambda)$ in which each vertex $u\in V = [n]$ is assigned a community $\sigma^\star_u\in [q]$ uniformly at random and then each possible edge $(u,v)$ is included independently with probability $\frac{a}{n}$ if $\sigma^\star_u = \sigma^\star_v$ and probability $\frac{b}{n}$ otherwise. Given a graph $G$ sampled from this distribution, the goal is to recover the underlying community membership $\sigma^\star\equiv (\sigma^\star_u)_{u\in V}$.

Throughout this paper, we allow the parameters $(a,b,q)\equiv (a_n,b_n,q_n)$ to depend on $n$, though we drop the subscripts for simplicity. In particular, we allow $q_n \equiv n^{\chi_n}\geq 2$ to grow with $n$ where $\chi\equiv \chi_n\in [0,1]$ controls the growth of the number of communities. We are most interested in the \emph{sparse regime}, where the average degree of the graph is constant. For this reason, we often reparameterize $(a,b)$ with 
\begin{equation}\label{eq: d,lambda}
d\equiv d_n=\frac{a+(q-1)b}{q}, \quad\quad \la\equiv \la_n =\frac{a-b}{a+(q-1)b}.
\end{equation}
Here, $d$ is the average degree of the block model and $\la$ is a measure of signal strength. With this parametrization, the relevant regime for the study of recovery is when $d, \la=\Theta(1)$. Note that when $q=q_n\to\infty$, $\lambda=\Theta(1)$ implies that $\lambda\geq 0$ (assortative) since $\lambda\in [-(q-1)^{-1},1]$ by definition. In our original parameters, this translates to $a=\Theta(q)$ and $b=\Theta(1)$. This will be the scaling assumed for the rest of the paper, unless otherwise noted. Recall that the Kesten--Stigum (KS) threshold is $d\lambda^2=1$. 

\begin{figure}
    \centering
    \includegraphics[width=0.325\linewidth]{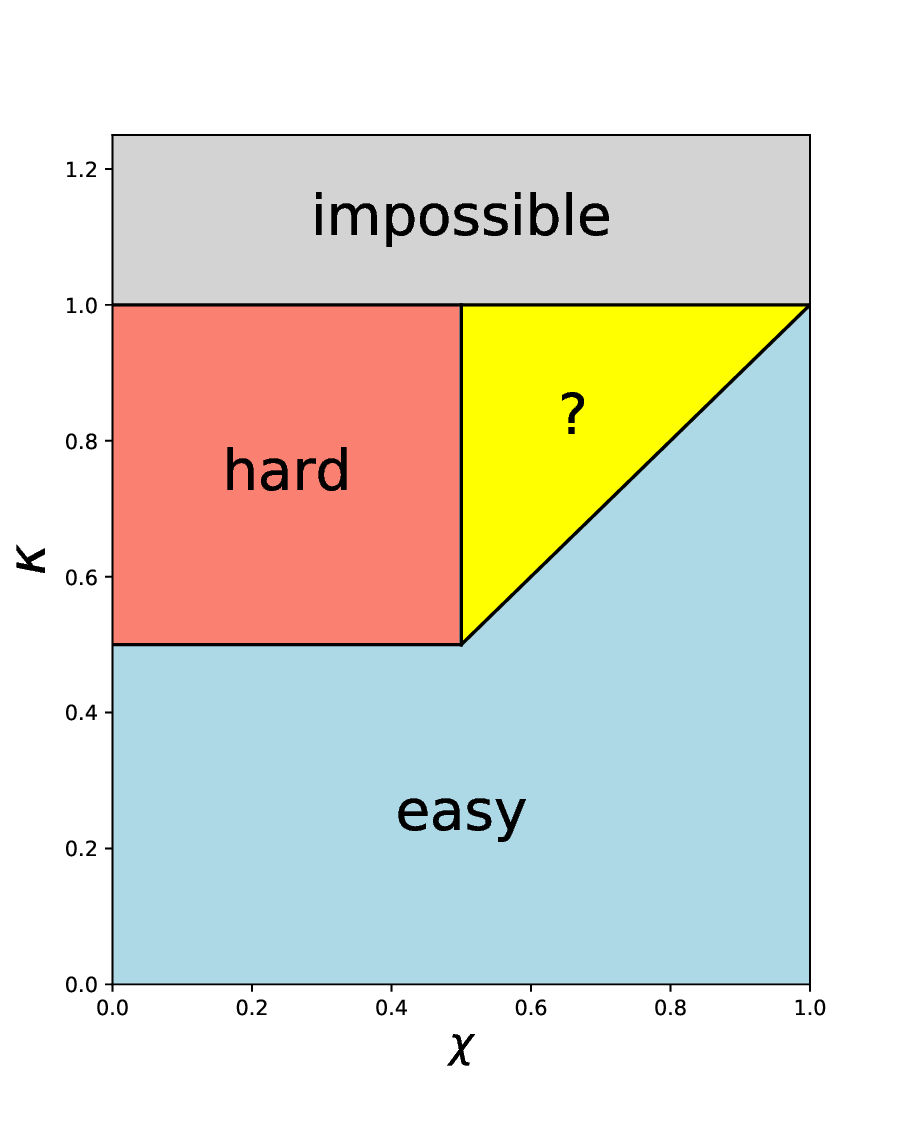}
    \includegraphics[width=0.325\linewidth]{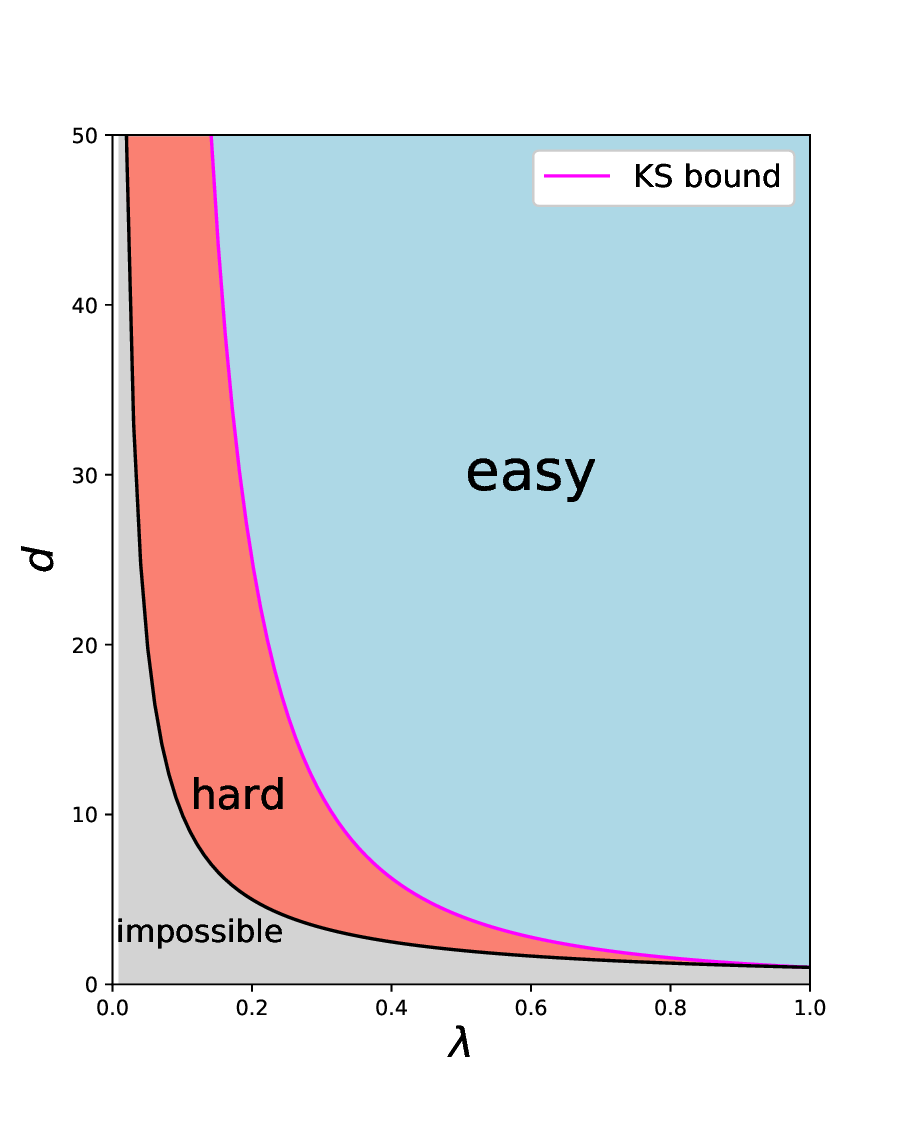}
    \includegraphics[width=0.325\linewidth]{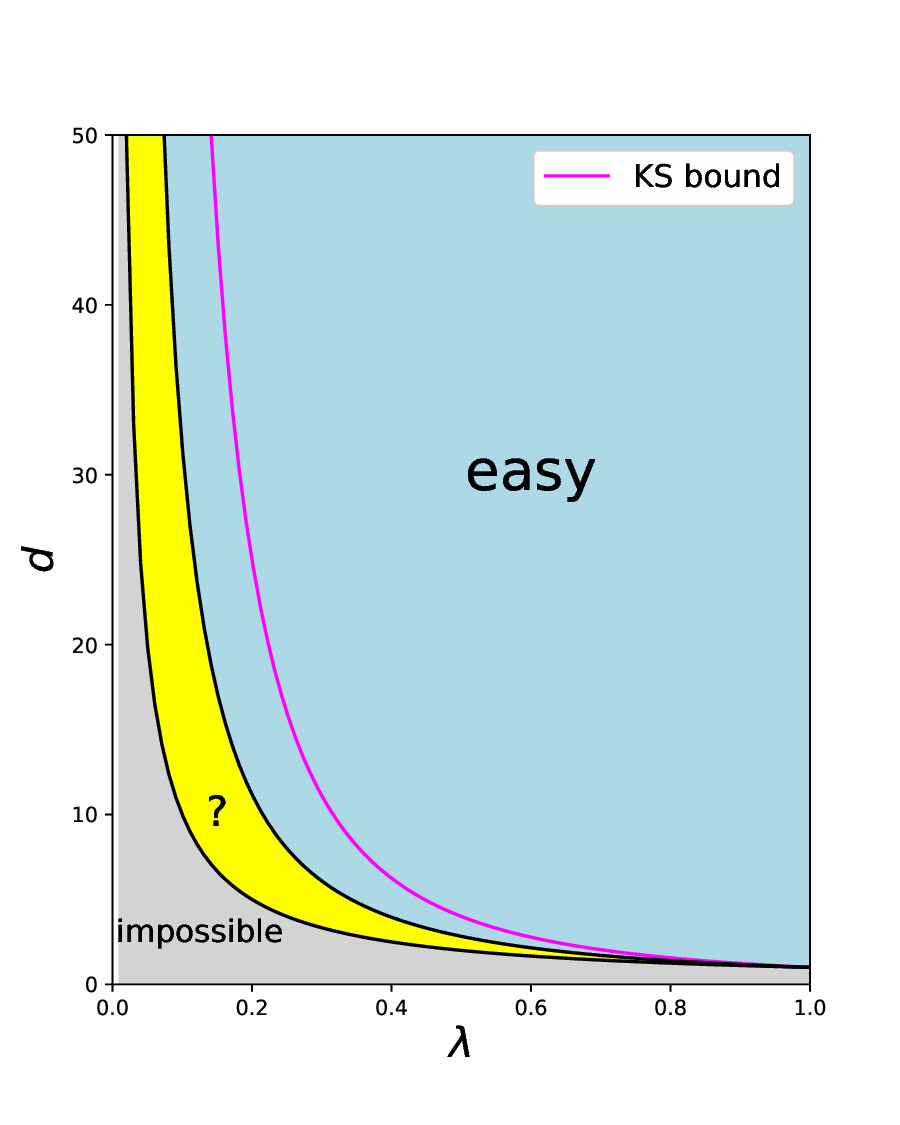}
    \caption{The blue regions represent computationally efficient regimes, the red regions represent low-degree hardness, the gray regions represent information-theoretic impossibility, and the yellow regions are still unknown. The left plot shows $\kappa$ vs $\chi$ where $q \asymp n^\chi$ and $\lambda \asymp d^{-\kappa}$. The middle and right plots show $d$ vs $\lambda$ for $q = n^{1/3}$ and $q = n^{2/3}$ respectively.}
    \label{fig:enter-label}
\end{figure}

Our results can be summarized by the diagrams in Figure~\ref{fig:enter-label}. In the following subsections we present the results in each region.

\subsubsection{Achievability of efficient recovery}
To measure the distance between two community assignments, we use a notion of \emph{alignment}, which is the maximum overlap of the two assignments after accounting for possible permutations of the labels. We will be interested in finding algorithms which can output a community assignment that has non-trivial alignment with the true assignment $\sigma^\star$. 
\begin{definition}\label{def: Align}
    The \textit{alignment} between two community assignments $\sigma,\tau \in [q]^V$ is defined as 
    \[ A(\sigma, \tau) \vcentcolon= \max_{\pi \in S_q} \frac{1}{n}\sum_{v=1}^n \mathbbm{1}\{ \sigma_v = \pi(\tau_v)\}. \]
    For $\delta > 0$ we say that an algorithm achieves \textit{$\delta$-recovery} if it outputs a partition $\widehat{\sigma} \in [q]^V$ such that $A(\widehat\sigma, \sigma^\star) > \delta$ with probability $1-o(1)$ as $n \to \infty$. We also refer to this notion as \textit{weak recovery} when no parameter $\delta$ is specified.
\end{definition}
The definition of $A(\sigma,\tau)$ has appeared in previous literature~\cite{CKPZ:18, AS:18,MSS:22, MSS24} where a fixed number of communities was studied. In that setting, for  $A(\widehat\sigma, \sigma^\star)$ to represent non-trivial weak recovery, it should exceed $1/q + \delta$ for a constant $\delta > 0$, as $1/q$ is the expected alignment achieved by a random assignment. However, because our work considers the case $q \to \infty$, the difference between $1/q + \delta$ and $\delta$ becomes negligible, and we use the simpler $\delta$ threshold.

\begin{theorem}\label{thm: main alg}
    Let $G \sim \SBM(n, q, d, \lambda)$. Then, the following holds.
    \begin{enumerate}
        \item Suppose $q=n^{\chi}$ where $1/2+\eps \leq \chi < 1$ for some $\eps>0$. Then, there exists a universal constant $C>0$ such that if $d\lambda^{1/\chi} > C(\log d)^2$ then Algorithm~\ref{alg: below KS} achieves $\delta$-recovery for some $\delta \equiv \delta(\eps,d)>0$.
        \item Suppose $q, \frac{n}{q\log^3 n} \to \infty$. If $d\la^2\geq 1+\eta$ for some $\eta>0$, then Algorithm~\ref{alg: above KS} achieves $\delta$-recovery for some $\delta\equiv \delta(\eta,d)>0$.
    \end{enumerate}
\end{theorem}
We will formally introduce the algorithms in Section~\ref{sec: Algorithmic} but remark here that they are based on computing non-backtracking walk statistics, and they run in polynomial time. In Subsection~\ref{subsec: discussion}, we will discuss our proof methods and compare our results to previous work in similar regimes, including~\cite{LR:15,GV:16,GMZZ:17,SZ:24,AS:15a,CX:16}.

Theorem~\ref{thm: main alg} shows a dichotomy when the number of communities is greater than $\sqrt{n}$ versus when the number of communities is less than $\sqrt{n}$. Below $\sqrt{n}$ communities we achieve $\delta$-recovery above the KS threshold, matching the behavior found in the constant number of communities regime. On the other hand, when there are more than $\sqrt{n}$ communities, we are able to achieve $\delta$-recovery below the KS threshold. Our result moreover indicates that a modified threshold based on the signal to noise ratio $d\lambda^{1/\chi}$ is a natural location for the phase transition in this regime. We ask a precise question in Subsection~\ref{subsec: discussion}.

We remark that in the case of fixed $q$, the result of~\cite{AS:15a} achieves weak recovery efficiently, above the KS bound. On the other hand, if $q = \Theta(n)$, weak recovery is trivial as the singleton communities comprise a constant fraction of the vertices. Thus, the conditions in part~2 represent the interesting regime of parameters.

\subsubsection{Low-degree hardness of recovery}
\label{subsec:low:degree}
Having established achievability results, a natural question arises: could the KS bound be crossed efficiently even for $q\ll n^{1/2}$ communities? We establish a hardness result suggesting that this is not possible. This is a question of \emph{average-case} complexity, where we aim to rule out an algorithm that solves \emph{most} random instances from a specific distribution, in contrast to the classical (worst-case) theory of NP-hardness where we aim to rule out an algorithm that solves \emph{all} instances.

We follow a popular paradigm in average-case complexity of statistics by studying \textit{low-degree polynomials}~\cite{HKPRSS:17,HS:17,H:18,KWB:19}. The original work on this framework considers hypothesis testing, and we consider the extension to \emph{recovery} problems, which originates from~\cite{SW:22}.

We represent the input graph $G \sim \SBM(n, q, d, \lambda)$ by the collection of $\binom{n}{2}$ binary variables $\{Y_{ij} \in \{0,1\} \,:\, 1 \le i < j \le n\}$, with $Y_{ij}$ encoding the presence (1) or absence (0) of edge $(i,j)$. For some bound $D \ge 0$, we consider estimators $f: \{0,1\}^{\binom{n}{2}} \to \R$ that are multivariate polynomials of degree at most $D$ in the variables $Y_{ij}$. This class of estimators is denoted $f\in \R[Y]_{\leq D}$. Consider the modest task of determining whether two distinct vertices (say, vertices 1 and 2 without loss of generality) belong to the same community or not. The success of degree-$D$ polynomials in such a task is measured by the \textit{degree-$D$ correlation}
\begin{equation}\label{eq:def:corr}
\Corr_{\le D} := \sup_{f\in \R[Y]_{\leq D}} \frac{\EE[f(Y) \cdot x]}{\sqrt{\EE[f(Y)^2]\cdot \EE[x^2]}}~~~\textnormal{where}~~~x := \one\{\sigma^\star_1=\sigma^\star_2\}-\frac{1}{q}.
\end{equation}
Here, all expectations are with respect to the joint distribution of $(x,Y)$ drawn from the block model $\SBM(n, q, d, \lambda)$. The normalization guarantees that $\Corr_{\leq D}\in [0,1]$, where $\Corr_{\le D} = 1-o(1)$ indicates near-perfect recovery of $x$, while $\Corr_{\le D} = \Omega(1)$ indicates non-trivial correlation with $x$, and $\Corr_{\le D} = o(1)$ indicates failure to correlate with $x$. It is a standard fact (see e.g.~\cite[Fact~1.1]{SW:22}) that $\Corr_{\le D}$ is directly related to the minimum mean squared error in estimating $\one\{\sigma^\star_1=\sigma^\star_2\}$ by a degree-$D$ polynomial.

Bounds on $\Corr_{\le D}$ have a conjectural interpretation for time complexity: if $\Corr_{\le D} = o(1)$ for some super-logarithmic degree, i.e.~$D \ge (\log n)^{1+\Omega(1)}$, this suggests that any algorithm needs super-polynomial time to obtain non-trivial correlation with $x$ (see~\cite[Conjecture~2.2.4]{H:18}). Furthermore, if $\Corr_{\le n^c} = o(1)$ for a constant $c > 0$, this suggests that runtime $\exp(n^{c-o(1)})$ is required (which is roughly the runtime needed to evaluate a degree-$n^c$ polynomial term-by-term). These rules are justified by their success in tracking the best known time complexity for a variety of different statistical models, along with the fact that degree-$O(\log n)$ polynomials tend to capture some of our best algorithmic tools like spectral methods. We note that there are limits to the style of models where the above conjectures hold up, but the SBM is considered a ``core'' problem in this area, where low-degree polynomials have already been ``stress-tested.'' In particular, for fixed $q$, the low-degree threshold coincides with the KS threshold for both hypothesis testing~\cite{HS:17,BBKMW:21} and weak recovery~\cite{SW:24,sbm-reduction}. The techniques of~\cite{SW:24} extend to $q \ll n^{1/8}$. A different approach relies on a low-degree hardness conjecture to establish average-case hardness of weak recovery for $q = n^{o(1)}$~\cite{sbm-reduction}. Coarser low-degree lower bounds for recovery are given by~\cite{LG23}, which apply for all $q \ll n$ but are not sharp enough to capture the precise KS transition.

Our hardness results show that if $q\ll n^{1/2}$ then polynomials of degree $n^{c}$ fail to achieve non-trivial correlation with $\one\{\sigma^\star_1=\sigma^\star_2\}$ below the KS bound. Notice that this theorem holds beyond the sparse regime, and even with constant edge probabilities with a modified SNR.
\begin{theorem}\label{thm: main hard}
Let $G \sim \SBM(n, q, a, b)$ for $a\leq \frac{n}{2}$. Suppose that there exist constants $\eps, \delta, \eta>0$ such that $q\leq n^{1/2-\eps}$, $a \geq \delta b$, and 
\[   \text{SNR}\equiv \frac{(a-b)^2}{q\left(a(1-a/n)+(q-1)b(1-b/n)\right)}\leq 1-\eta.   \]
Then, there exist constants $c\equiv c(\delta, \eta)>0$ and $C\equiv C(\eta)>0$ such that if $D\leq n^{c \, \eps}$ then
\begin{equation}\label{eq:thm:corr}
\Corr_ {\leq D}\leq C\sqrt{\frac{q}{n}}.
\end{equation}
\end{theorem}
By specializing our proof to the sparse SBM with growing number of communities, we establish low-degree hardness for all regions at or below the KS bound.

\begin{theorem}\label{thm:spase:mainhard}
    Consider the sparse SBM with growing number of communities where $q=\omega(1)$, $a=\Theta(q)$, and $b=\Theta(1)$. If $q=n^{1/2-\Omega(1)}$ and $d\la^2\leq 1$, then $\Corr_{\le n^{\Omega(1)}}\leq n^{-\Omega(1)}$.
\end{theorem}

Neither of the conditions $q=n^{1/2-\Omega(1)}$ and $d\la^2\leq 1$ are removable, as explained in the remark following Corollary~\ref{cor: sparse transition}. 

Observe that Theorem~\ref{thm: main hard} covers both the dense and the sparse regime. One motivation to consider the dense regime is to apply Theorem~\ref{thm: main hard} for $a,b\asymp n$ in order to establish computational barriers to \emph{graphon estimation}. This application is illustrated in Section~\ref{subsec:graphon}. Furthermore, observe that in this linear density regime $a,b\asymp n$, the condition $\text{SNR}\leq 1$ is stricter than $d\la^2\equiv \frac{(a-b)^2}{q(a+(q-1)b)}\leq 1$. We suspect that such correction is necessary simply because the complement of $G\sim \SBM(n, q, a, b)$ is distributed according to $\SBM(n,q,n-a,n-b)$ so that the ``correct'' signal-to-noise ratio must be invariant under the transformation $(a,b)\to (n-a,n-b)$. 

Unlike the other conditions in Theorem~\ref{thm: main hard}, the assumption $a\geq \delta b$ is likely an artifact of our proof. It would be interesting to remove it, as doing so could lead to important applications, such as in planted coloring where $a=0$.

\subsubsection{Phase transition at the KS bound}

Our algorithms give rise to low-degree estimators, so our results in the previous two subsections imply a sharp phase transition for low-degree polynomials at the KS threshold when $q \ll \sqrt{n}$. Combining our methods to prove Theorem \ref{thm: main alg}(2) with Theorem \ref{thm: main hard} yields the following result. 
\begin{corollary}\label{cor: sparse transition}
    Fix $\epsilon, \gamma, \eta > 0$. Let $G \sim \SBM(n, q, d, \lambda)$ where $q \to \infty$ and $q \leq n^{1/2-\epsilon}$. Then we have the following:
    \begin{enumerate}
        \item If $d\lambda^2 = 1 + \eta > 1$, then there exist constants $C \equiv C(\eta) > 0$ and $\delta \equiv \delta(\eta, d) > 0$ such that $\Corr_{\leq C\log n} \ge \delta$ for large enough $n$.
        \item If $d\lambda^2 \leq 1 - \eta$, then there exists a constant $c \equiv c(\eta) > 0$ such that $\Corr_{\leq n^{c \epsilon}} = o(1)$.
    \end{enumerate}
\end{corollary}
We remark that the proof methods for Theorem \ref{thm: main alg}(1) imply that this transition breaks down for $q \gg \sqrt{n}$, where $\Corr_{\leq C\log n} \ge \delta$ even when $d\lambda^2 < 1$ but $d\lambda^{1/\chi} > C(\log d)^2$ (see Section~\ref{subsec: LDI}).

The low-degree threshold for $q \gg \sqrt{n}$ remains an open question, and we do not give a low-degree lower bound in this regime. The existing lower bound of~\cite{LG23} covers the full range of $q \ll n$, and when specialized to our setting where $b,d=\Theta(1)$ and $a,q \to \infty$, this shows low-degree hardness of recovery when $\lambda^2 \ll \min\{1,n/q^2\}$ where $\ll$ hides $\log(n)$ factors. For $q \ll \sqrt{n}$ this matches the KS threshold up to log factors, and for $q \gg \sqrt{n}$ this gives a lower bound that is non-trivial but far from our positive result $d\lambda^{1/\chi} > C(\log d)^2$ since $\lambda$ needs to shrink with $n$.

\subsubsection{Computational barriers to graphon estimation}
\label{subsec:graphon}

The stochastic block model is a specific example of a \textit{graphon}~\cite{L:12}, so our hardness result in Theorem~\ref{thm: main hard} has a natural implication in graphon estimation. Specifically, consider a graph with $n$ vertices where the edges are generated independently connecting vertices $i\neq j$ with probability 
\[
M_{i,j}= f(\xi_i,\xi_j).
\]
Here, $(\xi_1,\ldots, \xi_n)$ are i.i.d.\ latent feature vectors sampled from an unknown distribution $\P_{\xi}$ supported on $[0,1]$. The function $f:[0,1]\times [0,1]\to \R$ is a symmetric function called a graphon. We denote the adjacency matrix of the resulting observed graph by $Y\in \{0,1\}^{n\times n}$ where $Y_{i,i}=0$ by convention. In the case where $\P_{\xi}$ has discrete support, or equivalently where $f$ is blockwise constant, this observed graph specializes to the (possibly non-symmetric) stochastic block model. 

\textit{Graphon estimation} has gained great interest in the last decade \cite{OW:13, ACC:13, AC:14, GLZ:15, KTV:17}. The goal is to estimate the graphon $f$ given the adjacency matrix $Y$. In particular, it was shown in \cite[Section~3]{KTZ17} that estimating $f$ in squared error can be achieved by estimating $M$ in Frobenius norm:
\[
\ell(\widehat{M},M)=\binom{n}{2}^{-1}\sum_{1\leq i<j\leq n} (\wM_{i,j}-M_{i,j})^2,
\]
where $\wM \in [0,1]^{n\times n}$. Here, we focus on the SBM class with $q$ blocks
\[
\cM_q=\Big\{M\in [0,1]^{n\times n}: M_{i,j}=Q_{z_i,z_j}\one\{i\neq j\}~~\textnormal{for some}~~Q=Q^\top \in [0,1]^{q\times q},\, z\in [q]^n\Big\}.
\]
The best known polynomial time estimator $\wM$ for graphon estimation in the SBM class is \emph{universal singular value thresholding (UVST)} by~\cite{Cha15}. The works of~\cite{Xu18} and~\cite{KV19} independently established the error bound
\[
\sup_{M\in \cM_q} \EE\Big[\ell(\wM_{\textnormal{UVST}},M)\Big]\leq C\frac{q}{n},
\]
for some universal constant $C>0$. In contrast, it was shown in \cite{GLZ15} that the min-max error rate, ignoring computational constraints, is
\[ \inf_{\wM}\sup_{M\in \cM_q} \EE\big[\ell(\wM_{\textnormal{UVST}},M)\big]\asymp \frac{q^2}{n^2}+\frac{\log q}{n}, \]
which points to a computational-statistical gap. Recently, rigorous evidence towards such a gap was obtained by~\cite{LG23}, where they showed that if $2\leq q \leq \sqrt{n}$ then there is a universal constant $c$ such that 
\[
\inf_{\wM\in \R[Y]_{\leq D}^{n\times n}}\sup_{M\in \cM_q} \EE\Big[\ell(\wM,M)\Big]\geq \frac{cq}{nD^4},
\]
where $\wM\in \R[Y]_{\leq D}^{n\times n}$ means $\wM_{i,j}\in \R[Y]_{\leq D}$ for each $i,j \in [n]$. Although this result achieves the rate $\frac{q}{n}$ up to log factors when $D=\log^{1+\eps}n$, it is not clear if the same rate can be achieved for $D=n^{\delta}$ for some fixed constant $\delta>0$ nor that there is a tradeoff between the degree (or equivalently runtime) and the estimation error due to the $D^4$ factor. Our corollary removes this $D^4$ factor when $q\leq n^{1/2-\eps}$.
\begin{corollary}\label{cor: graphon estimation}
    Suppose $2\leq q\leq n^{1/2-\eps}$ for some $\eps>0$. There exists a constant $c \equiv c(\eps)$ and a universal constant $C>0$ such that if $D\leq n^{c}$ then
    \[
    \sup_{M\in \cM_q}\inf_{\wM\in \R[Y]_{\leq D}^{n\times n}} \EE\Big[\ell(\wM,M)\Big]\geq \frac{Cq}{n}.
    \]
\end{corollary}

\noindent As discussed in Section~3.1 of~\cite{LG23}, when $q\geq n^{1/2+\eps}$ the symmetric block model, i.e., when the connection probabilities are the same within communities and across communities, fails to imply sharp hardness for graphon estimation. To obtain sharp hardness, it appears necessary to construct a different stochastic block model that is harder than the symmetric block model, and this remains an open problem.

\subsubsection{Information-theoretic thresholds and detection}
\label{subsec:info}
To complement the algorithmic picture we study the information-theoretic behavior as well. We show that the scaling for the threshold is determined by the SNR $d\lambda$ whenever the number of communities increases with the size of the graph. 
\begin{theorem}\label{thm: IT}
    Let $G \sim \SBM(n,q,d,\lambda)$ with $q = \omega(1)$. Then the following holds:
    \begin{enumerate}
        \item Suppose $d\lambda > 14$. Then the exponential-time~Algorithm~\ref{alg: inefficient} achieves $\eps$-recovery, where $\eps=\Omega(1)$ as $n\to\infty$. 
        \item Suppose $d\lambda < 1$ and $q = o\left(n\right)$. Then it is information-theoretically impossible to to achieve $\delta$-recovery for any fixed $\delta > 0$.
    \end{enumerate}
\end{theorem}

Finally, we show that the question of detection is easy in this regime, as a simple counting algorithm succeeds well below the recovery threshold.
\begin{theorem}\label{thm: detection}
    Let $q = \omega(1)$ and fix $d, \lambda \neq 0$. Then counting triangles distinguishes between $G \sim \SBM(n,q,d,\lambda)$ and $G(n, \frac{d}{n})$ with probability $1-o(1)$.
\end{theorem}

\subsection{Proof techniques} 
Our proof for the algorithmic feasibility of weak recovery is centered on a weighted non-backtracking walk statistic. The algorithm is inspired by the non-backtracking walk counting argument of~\cite{MNS:18}. The proof of correctness has two main steps: first we show that the non-backtracking walks give a good estimate for whether or not two fixed vertices are in the same community or not. This is similar to a strategy used above the KS bound when $q = 2$, but surprisingly we are able to implement this step below the KS bound as well when $q\gg\sqrt{n}$. We provide more details of this unexpected step below. In the second step, we assemble these pairwise estimates into a global partition of the vertices using a coupling with trees.

The first step is to show that the non-backtracking walk statistics concentrate well around their means. We approach this using the second moment method. Our proof combines combinatorial bounds on the number of non-backtracking walks with a moment analysis of the weights of these walks. The concentration allows us to reduce the analysis of the complicated non-backtracking walk statistics to the analysis of a local estimator depending only on the community labels in the neighborhoods of the starting and ending vertices. 

In the regime above the KS threshold, a careful second moment computation is strong enough for this reduction to follow from a straightforward union bound. Surprisingly, and in contrast to the case of fixed $q$, we are able to use the second moment effectively below the KS threshold as well. For fixed $q$ below the KS bound the second moment bound is far too weak to deduce any useful concentration. While carrying out the analysis we discover that it is possible to exploit the size of $q$ to obtain a moderate concentration. Essentially, the growth of $q$ can compensate for the terms that typically make the variance large below the KS threshold. By itself this moderate concentration suffices to prove $n^{-c}$-recovery for some $\frac{1}{q} \ll n^{-c} \ll 1$. This is already an unexpected new behavior that arises due to the growth of $q$. In order to achieve our desired $\delta$-recovery for some fixed $\delta > 0$, we bootstrap this moderate concentration by partitioning the graph into many pieces. This effectively gives us independent trials of moderate concentration, which when combined boost the concentration to the same quantitative strength as above the KS bound. This is also the source of the $(\log d)^2$ term in our bound, as partitioning the graph reduces the SNR which we require to be large enough on each piece. 

We next estimate the performance of the local estimators, which are effectively 
a weighted count of how many pairs of vertices in the neighborhoods have matching labels. The analysis of the local estimators typically relates to the theory of broadcast processes on Galton--Watson trees, but turns out to be simpler in our setting, where it suffices to restrict our attention to events that depend only on the depth-1 neighborhoods of vertices. This shows that the local estimators, and thus also the non-backtracking walk statistics, correlate with the indicator that start and end vertices have matching labels.

Finally, the signal from the pairwise estimators is assembled into a global partition. This is not an obvious step, as for example, the analysis of the pairwise estimators does not rule out the following situation: with probability $\delta$ we recover the communities exactly and otherwise we recover nothing. This would yield an algorithm that outputs a partition with $\delta$ alignment with the true communities in expectation, but not with high probability. We circumvent this via a standard coupling of neighborhoods in a sparse random graph with an independent collection of depth-1 random trees. The independence here is crucial, as it precludes strong dependencies between the local estimators that are a priori arbitrarily correlated.

In order to establish the low-degree hardness result of Theorem~\ref{thm: main hard}, we build on the recent proof framework of~\cite{SW:24}, which analyzes the SBM with a fixed number of communities $q=O(1)$. However, a direct application of their methods only handles the regime $q\ll n^{1/8}$. Thus, the central technical challenge is to extend these techniques and prove hardness below the KS bound when $q\ll n^{1/2}$. Roughly speaking, we overcome this hurdle by choosing a different basis for $\R[Y]_{\leq D}$ to expand $\EE[f(Y)\cdot x]$ in Eq.~\eqref{eq:def:corr}. This new basis enables explicit calculations and, in turn, allows us to cover the full regime $q\ll n^{1/2}$. We refer to Section~\ref{sec: LDH} for more details.

We briefly introduce the ideas for our information-theoretic bounds as well. The achievability result in Theorem \ref{thm: IT} begins with a brute-force search over all partitions using a strategy of Banks, Moore, Neeman, and Netrapalli~\cite{BMNN:16} to obtain a polynomially small weak recovery. Then, we perform one pass of depth-1 belief propagation to boost this to recover a constant fraction of the labels. Conversely, the impossibility result follows from a percolation-style argument claiming that the ``propagation'' of signal from a vertex dies out quickly.

\subsection{Discussion and future directions}\label{subsec: discussion}

In this subsection we discuss how our Theorem~\ref{thm: main alg} relates to prior literature. Much of the literature deals either with the setting of slowly growing $q$, or is interested in the denser regime where the average degree is logarithmic. We nonetheless introduce the relevant settings and compare our results and proof techniques in the overlapping regimes. 

Lei and Rinaldo~\cite{LR:15} use spectral clustering to achieve weak recovery. Their results are about the semi-sparse regime where $a \geq \log n$, and their error rate is meaningful when $d\lambda^2 > C(1+\lambda(q-1))$. This matches our result up to constant factors when $\lambda q$ is rather small, corresponding to the case where $q \ll \sqrt{d}$. This requires the degree to be large when there are many communities. 

Another spectral method using the Bethe Hessian matrix, initially introduced in \cite{SKZ:14}, has seen an increase in popularity for performing weak recovery. The recent work of Stephan and Zhu~\cite{SZ:24} uses the spectral properties of this matrix to achieve weak recovery for $d\lambda^2 > 1+\epsilon$ in the regime where $d,q \leq \mathrm{polylog}(n)$. This matches our result for this regime of constant or slowly growing $d$ and $q$. They also give sufficient conditions to achieve $(1-o(1))$-recovery, in which they require $d \gg q$ and $d\lambda^2 \gg q$. 

On the other hand Gu\'edon and Vershynin~\cite{GV:16} use semidefinite programming methods to study community detection in sparse (constant average degree) stochastic block models. In our setting, their results can handle all $q$ and translate to a requirement that $d\lambda^2 > C$ where $C$ is a large constant that has an explicit dependence on the error rate of the algorithm. However, the guarantee of their algorithm is a matrix with small $L^2$ distance, and it is not clear how to extract a community partition from this guarantee.

We mention that the original work of Abbe and Sandon~\cite{AS:15a} achieving the KS threshold for all fixed $q$ uses a novel Sphere Comparison algorithm and has a quantitative error rate of $C \exp(-d\lambda^2/q^3)$. In our setting this requires the much stronger condition $d\lambda^2 = \Omega(q^3)$ to be successful.

Gao, Ma, Zhang, and Zhou~\cite{GMZZ:17} study a slightly more refined question of recovery --- the question of achieving the optimal fraction of community labels --- also using spectral methods. Naturally their conditions are stronger than ours, requiring that $d\lambda^2 > Cq(1+\lambda(q-1))$ at initialization. 

We also mention here that Chen and Xu~\cite{CX:16} studied the question of exact recovery in a more general setting. They are able to derive information-theoretic impossibility in a regime of the form $d\lambda^2 \leq \frac{\log n}{q}$ or $a \leq q\log n$. On the other hand, their efficient algorithm succeeds when $d\lambda^2 \geq C(\frac{1+\lambda(q-1)}{q} \cdot \log n + 1)$. When $q \gg \log n$ and $\lambda \ll \frac{1}{\log n}$ this matches our threshold up to constant factors.

Our results also leave many natural directions for future work in this new regime. We believe that the SNR $d\lambda^{1/\chi}$ is the right quantity to determine the feasibility of efficient recovery.
\begin{question}
Let $q=n^{\chi}$, where $\chi>\frac{1}{2}+\Omega(1)$. Does the phase transition for efficient weak recovery occur at $\lambda = \widetilde{\Theta}(d^{-\chi})$, where $\widetilde{\Theta}$ hides poly-log factors? Moreover, is there a sharp or coarse transition in the recovery accuracy as $\la$ varies as a function of $d$?
\end{question}
We refer to sharp versus coarse thresholds in the sense of Friedgut~\cite{F:99}. Our work answers the achievability side of this question, showing that $\lambda > Cd^{-\chi}(\log d)^{2\chi}$ suffices for efficient weak recovery, and it would be interesting to provide a matching lower bound for efficient weak recovery within the low-degree framework.

It is also natural to ask if similar statements are true for denser block models with $q \gg \sqrt{n}$ communities. This conjecture corresponds to the yellow regions in the leftmost diagram in Figure~\ref{fig:enter-label} being low-degree hard. Conversely, handling the logarithmic factors in the achievability result seems to require a new idea. Our analysis leaves open the ``critical'' regime where $q \asymp \sqrt{n}$. It would be interesting to study the question of weak recovery in this critical window as well.

We note that some recent work~\cite{gmm} on clustering high-dimensional Gaussian mixture models has some parallels to ours. This problem shares some similarities to the SBM, where the number of communities $q$ corresponds to the number of clusters, and the number of vertices $n$ corresponds to the number of samples. The work~\cite{gmm} identifies the computational threshold at an analogous KS bound up to log factors, giving a low-degree lower bound that matches the known poly-time algorithms. Curiously, the behavior changes once $q$ exceeds $\sqrt{n}$, just like in the SBM. When $q \gg \sqrt{n}$, the optimal poly-time algorithm for clustering is a simple one based on grouping nearby points~\cite[Appendix~F]{gmm}. A non-rigorous transposition (supported by e.g.,~\cite{GV:19}) suggests a similar change in computational threshold for rather dense block models. However, this transposition breaks down for constant $d$ and it is not clear that a similarly simple algorithm works. Since we focus on the sparse setting of the SBM we leave this potential connection to future research.

\subsection{An intuition from the dense setting which gives the SNR \texorpdfstring{$d \lambda^{1/\chi}$}{dl\^1/chi}}\label{subsec: dense heuristic}
To motivate the new SNR, we first consider the performance of a simple algorithm based on counting common neighbors in the dense case and see how this leads to the quantity $d\lambda^{1/\chi}$. For this intuition, we assume a setup of parameters in which $\frac{n}{q}$ is an integer and all communities have equal sizes.

Fix two vertices $u, v$, and let $X$ be the number of common neighbors of $u$ and $v$. The distribution of $X$ can be described as follows:
\[ X \sim \begin{cases}
    \mathrm{Bin}\left(\frac{n}{q}, \frac{a^2}{n^2}\right) + \mathrm{Bin}\left(n\left(1 - \frac{1}{q}\right), \frac{b^2}{n^2}\right) & \sigma^\star_u = \sigma^\star_v, \\
    \mathrm{Bin}\left(\frac{2n}{q}, \frac{ab}{n^2}\right) + \mathrm{Bin}\left(n\left(1 - \frac{2}{q}\right), \frac{b^2}{n^2}\right) & \sigma^\star_u \neq  \sigma^\star_v.
\end{cases} \]
In the above, the binomial random variables are assumed to be independent. A straightforward calculation shows that the difference in expectations is $\frac{(a-b)^2}{nq}$. Applying the Chernoff bound for binomial random variables, we have the following tail estimates
\begin{align*}
    \Prob{\mathrm{Bin}\left(\frac{2n}{q}, \frac{ab}{n^2}\right) > \frac{2ab}{nq} + \frac{(a-b)^2}{2nq}} &< \exp(-\frac{(a-b)^4}{nq(16ab + 2(a-b)^2)}) \\
    \Prob{\mathrm{Bin}\left(n\left(1 - \frac{2}{q}\right), \frac{b^2}{n^2}\right) > \frac{b^2(1-\frac{2}{q})}{n} + \frac{(a-b)^2}{2nq}} &< \exp(-\frac{(a-b)^4}{nq(8b^2(q-2) + 2(a-b)^2)}) \\
    \Prob{\mathrm{Bin}\left(\frac{n}{q}, \frac{a^2}{n^2}\right) < \frac{a^2}{nq} - \frac{(a-b)^2}{2nq}} &< \exp(-\frac{(a-b)^4}{8nqa^2}).
\end{align*}
For the final binomial, notice that its expectation $\frac{b^2(1-\frac{1}{q})}{n}$ is of smaller order than $\frac{(a-b)^2}{nq}$ as long as $q\lambda \to \infty$, so the associated lower tail event has probability 0. Thus, to union bound over $q-1$ copies of the first tail event and one copy of the third tail event, we need to satisfy the conditions 
\[ \frac{d^2\lambda^2q}{n} \gg \log n \quad \text{ and } \quad \frac{d^2q^2\lambda^3}{n} \gg \log n \quad \text{ and } \quad \frac{d^2q^2\lambda^4}{n} \gg \log n. \]
These conditions have a few enlightening interpretations. Notice that if $d\lambda^2 < 1$, then a necessary condition for this test to succeed is that $q \gg \sqrt{n}$. This complements the regime in which the low-degree hardness argument succeeds, and shows that the condition is natural.

Fascinatingly, this simple test also motivates the form of the modified signal-to-noise ratio. Suppose we assume the scaling $q = n^{\chi}$ for $\chi > \frac{1}{2}$ as in the sparse setting and that $d\lambda^c = \Theta(1)$. Under this assumption, we have $\lambda \asymp d^{1/c}$ and the conditions translate to
\[ d \gg n^{\frac{1-\chi}{2-2/c}} \quad \text{ and } \quad d \ll n^{\frac{2\chi-1}{2-4/c}}. \]
This interval is nonempty exactly when $c \geq \frac{1}{\chi}$. Heuristically, this suggests that the scaling $d\lambda^{1/\chi} = \Theta(1)$ is a natural scaling to consider. As we will see in the next subsection, this same quantity will be crucial in our bound on the regime for weak recovery for sparse block models.

\subsection{The SNR $d\lambda^{1/\chi}$ in the sparse setting}
We highlight how our new SNR arises from the second moment computation. We need a few definitions from Section \ref{sec: Algorithmic}. For any edge $e$ define its weight $W_e$ to be $\mathbbm{1}\{e \in G\} - \frac{d}{n}$. This extends to a weight for any path $\gamma$ defined as $X_\gamma = \prod_{e \in \gamma} W_e$. We will also use $s = d\lambda$. The following is a simplified combination of Propositions~\ref{prop: path count moments} and~\ref{prop: path count moments KS}. 
\begin{proposition}
    For any two fixed vertices $u,v$, define the \textit{weighted self-avoiding walk count} between $u$ and $v$ as $S_{u,v} = \sum_{\gamma} X_\gamma$ where $\gamma$ sums over all self-avoiding walks of length $k = \lfloor \beta(d, s)\log n \rfloor$ from $u$ to $v$. Let $u, v, u', v'$ be distinct vertices. Then the following holds:
    \begin{align*}
        \E{S_{u,v}} &= \begin{cases}
            (1 + o(1)) \frac{(q-1) s^k}{n} & \sigma^\star_u = \sigma^\star_v, \\
            (1 + o(1)) \frac{-s^k}{n} & \sigma^\star_u \neq \sigma^\star_v,
        \end{cases} \\
        \E{S_{u,v}^2} &= \begin{cases}
            (1 + o(1)) \frac{s^2}{(s-1)^2} \cdot \frac{q^2 s^{2k}}{n^2} & \sigma^\star_u = \sigma^\star_v, \\
            (1+o(1))(\sum\limits_{\substack{i,j \in \mathbb{N}, i+j \leq k}} \frac{qs^{2k}}{n^2}\left(\frac{d}{s^2}\right)^{i+j} + \frac{d^k}{n}) & \sigma^\star_u \neq \sigma^\star_v.
        \end{cases}
    \end{align*}
\end{proposition}
In the above expression, the two dominant terms in the second moment when $\sigma^\star_u\neq \sigma^\star_v$, namely $\sum_{i+j\leq k} \frac{qs^{2k}}{n^2}\left(\frac{d}{s^2}\right)^{i+j}$ and $\frac{d^k}{n}$, arise from pairs of nearly disjoint walks and pairs of identical walks respectively. We refer to the proof of Proposition~\ref{prop: path count moments} for precise definitions of these pairs of walks. From this we can see where our phase transitions come from. Above the KS bound, or equivalently $\frac{d}{s^2} < 1$, the second moment is dominated by the term 
\[ \sum_{i+j\leq k} \frac{qs^{2k}}{n^2}\left(\frac{d}{s^2}\right)^{i+j} > \frac{s^2}{s^2-d} \cdot \frac{qs^{2k}}{n^2} \gg \frac{d^k}{n} \]
for a suitable choice of $k$. This is sufficient to carry out the concentration argument. On the other hand, below the KS bound, the second moment is dominated by 
\[ \frac{d^k}{n} \gg \frac{qd^k}{n^{2-o(1)}} \geq \sum_{i+j\leq k} \frac{qs^{2k}}{n^2}\left(\frac{d}{s^2}\right)^{i+j}. \]
In order to obtain our desired concentration, we need $\E{S_{u,v}^2} = o\left(\frac{q^2s^{2k}}{n^2}\right)$ or equivalently 
\[ \frac{d^k}{s^{2k}} \cdot \frac{n}{q^2} \ll 1. \]
We see here that since $\frac{d}{s^2} > 1$, this is only possible when $q \gg \sqrt{n}$ --- this is the reason for the transition in behavior at $\sqrt{n}$ communities. More precisely, this condition translates to $k < \frac{2\chi-1}{\log d/s^2} \log n.$
Notice that in order to have concentration of the random variable, we also require $\frac{qs^k}{n} \to \infty$, which translates to $k > \frac{1-\chi}{\log s} \log n$. Ensuring that these two restrictions on $k$ are non-contradictory is equivalent to 
\[ \frac{1-\chi}{\log s} < \frac{2\chi-1}{\log d/s^2} \iff (1-\chi)\log d < \log s \iff -\chi \log d < \log \lambda \iff d\lambda^{1/\chi} > 1 \]
which is the new SNR.

\subsection*{Organization}
In Section~\ref{sec: Algorithmic} we prove our algorithmic results --- Theorem~\ref{thm: main alg} and Corollary~\ref{cor: sparse transition}(1). Section~\ref{sec: LDH} is devoted to the hardness result Theorem~\ref{thm: main hard}, and Corollary~\ref{cor: graphon estimation} about graphon estimation. Finally, we conclude with the proof of our information-theoretic result, Theorem~\ref{thm: IT}, in Section~\ref{sec: IT} and our detection result, Theorem~\ref{thm: detection}, in Section~\ref{sec: detection}.

\subsection*{Acknowledgments}
EM and YS thank Allan Sly for helpful discussions. We thank the anonymous referees for feedback that improved this paper. BC is supported by an NSF Graduate Research Fellowship and Simons Investigator award (622132). EM was partially supported by NSF DMS-2031883, Bush Faculty Fellowship ONR-N00014-20-1-2826, Simons Investigator award (622132), and MURI N000142412742. ASW was partially supported by a Sloan Research Fellowship and NSF CAREER Award CCF-2338091.

\section{Achievability of efficient recovery}\label{sec: Algorithmic}

In this section we will prove our results concerning when weak recovery is efficiently possible. The upshot of this section is that weak recovery is efficiently possible above the KS threshold as in the finite community case, but remains possible below the KS threshold when $q$ is sufficiently large. The basic structure of the section is as follows. Subsection~\ref{sec: Alg} begins the analysis of the sparse setting and introduces the basic components of our algorithm, and gives a polynomial bound on the runtime for computing the relevant statistics. Subsection~\ref{sec: PC} recalls combinatorial path bounds that we use to prove our second moment bounds, which is carried out in Subsection~\ref{sec: WPS}. The remainder of the section introduces Algorithms~\ref{alg: below KS} and~\ref{alg: above KS} in detail and is devoted to combining the ingredients to prove Theorem~\ref{thm: main alg} and Corollary~\ref{cor: sparse transition}(1). 

\subsection{Algorithm preliminaries}\label{sec: Alg}
Our efficient algorithms will all follow the same structure. The specific algorithm used in each regime will differ slightly, so we postpone more detailed descriptions to the later subsections. The basic outline of our algorithm proceeds as follows: 
\begin{enumerate}[label=\arabic*.]
    \item Choose a set of vertices $U \subseteq V$ uniformly at random with size $O(n^{\frac{\chi+1}{2}})$. Collect the $q$ vertices of highest degree in $U$ to be $U_*$. These vertices will serve as representatives for the $q$ communities.
    \item For each $w \in V \setminus U$ and each $u_* \in U_*$, use the weighted non-backtracking walk count between $w$ and $u_*$ to decide whether or not the community represented by $u_*$ is a candidate for $w$.
    \item Choose one candidate representative for each $w$, and output the partition defined by the set of vertices assigned to the same representative.
\end{enumerate}

We analyze the initialization step (Step~1) of the outline here. 
\begin{lemma}\label{lemma: init}
    Suppose $q = o(\frac{n}{(\log \log n)^{\log\log n}})$. Let $U$ be a set of size $\sqrt{nq}$ and let $U_*$ be a uniformly random subset of size $q$ drawn from all vertices in $U$ which have degree at least $\log\log n$ to $V \setminus U$. With high probability, there are at least $\frac{q}{4}$ communities with a unique representative in $U_*$. 
\end{lemma}
\begin{proof}
    We reveal $\sigma^\star$ so that all edges in the graph are drawn independently from each other. Let $U_i$ be the vertices in $U$ with label $i$, and $V_i$ be the vertices in $V \setminus U$ with label $i$. With high probability, 
    \[ |U_i| = \sqrt{\frac{n}{q}} + O\left(\left(\frac{n}{q}\right)^{1/4}\log(\frac{n}{q})\right)\text{ and } |V_i| = \frac{n}{q} + O\left(\sqrt{\frac{n}{q}}\log(\frac{n}{q})\right) \quad \forall i. \]
    The degrees of vertices in $U$ to $V \setminus U$ are independent, as we have revealed all of the community labels. Moreover, each vertex $u_i \in U_i$ has degree distribution 
    \[ D_i = \mathrm{Bin}\left(|V_i|, \frac{a}{n}\right) + \mathrm{Bin}\left(n - \sqrt{nq} - |V_i|, \frac{b}{n}\right) \]
    where $+$ denotes convolution in this context. This differs from 
    \[ \widetilde{D_i} = \mathrm{Pois}\left(|V_i|\cdot \frac{a}{n} + \left(n - \sqrt{nq} - |V_i|\right)\cdot \frac{b}{n}\right) = \mathrm{Pois}\left(d + O\left(\sqrt{\frac{q}{n}}\log(\frac{n}{q})\right)\right) \]
    by $O\left(\frac{q}{n}\right)$ in total variation distance (see e.g. Lemma 5 in \cite{MNS:15}). We compare $\Prob{\widetilde{D_i} > \log\log n}$ with $\Prob{D > \log\log n}$ where $D \sim \mathrm{Pois}(d)$. First notice that 
    \[ \Prob{\widetilde D_i > \log\log n}\wedge \Prob{D > \log\log n} > \frac{1}{(\log\log n)^{\log\log n}} \]
    A standard Chernoff bound yields 
    \[ \Prob{\widetilde{D_i} > \log n} \vee\Prob{D > \log n} < \frac{(ed)^{\log n}e^{-d}}{(\log n)^{\log n}} = o\left(\frac{1}{(\log \log n)^{\log\log n}}\right). \]
    Thus we have 
    \[ \Prob{D > \log\log n} = (1+o(1))\sum_{k=\log\log n}^{\log n} e^{-d} \frac{d^k}{k!} \] 
    and 
    \[ \Prob{\widetilde D_i > \log\log n} = (1+o(1))\sum_{k=\log\log n}^{\log n} e^{-d-O(\sqrt{q/n}\log(n/q))} \frac{(d+O(\sqrt{q/n}\log(n/q)))^k}{k!}. \]
    For each $\log\log n < k < \log n$, we have 
    \begin{align*}
        \frac{e^{-d}\frac{d^k}{k!}}{e^{-d-O(\sqrt{q/n}\log(n/q))} \frac{(d+O(\sqrt{q/n}\log(n/q)))^k}{k!}} &= e^{O(\sqrt{q/n}\log(n/q))}\left(1 + O\left(\sqrt{\frac{q}{n}}\log(n/q)\right)\right)^k \\
        &\leq \exp(O\left(k\sqrt{\frac{q}{n}}\log(\frac{n}{q})\right)) = 1+o(1).
    \end{align*}
    Thus, all together we have shown that $\Prob{\widetilde D_i > \log\log n} = (1+o(1))\Prob{D > \log\log n}$ for every $i$. Now, the number of vertices in $U$ with label $i$ and degree larger than $\log\log n$ is distributed as 
    \[ \mathrm{Bin}\left(|U_i|, \Prob{D_i > \log\log n}\right). \]
    Thus, with high probability, there are 
    \[ |U_i|\cdot \Prob{D_i > \log\log n} = (1+o(1))\sqrt{\frac{n}{q}} \cdot \Prob{D > \log\log n} \]
    such vertices. Here we used that the expectation tends to infinity. Notice that the right hand side is independent from $i$. Thus, a randomly selected vertex of $U$ with degree at least $\log\log n$ has label $i$ with probability $\frac{1+o(1)}{q}$ for every $i$. If we let $X$ denote the number of communities with a unique representative selected in $U_*$, then 
    \[ \E{X} = q^2\cdot \frac{1+o(1)}{q} \cdot \left(1 - \frac{1+o(1)}{q}\right)^{q-1} = (1+o(1))\frac{q}{e} \]
    and 
    \[ \Var{X} = (1+o(1))\frac{q}{e} + q^2(q-1)^2\left(\frac{1+o(1)}{q}\right)^2\left(1 - \frac{1+o(1)}{q}\right)^{q-2} - (1+o(1))\frac{q^2}{e^2} = o(q^2). \]
    By Chebyshev's inequality we deduce that $X \geq \frac{q}{4}$ with high probability.
\end{proof}

We will restrict our focus to recovering those vertices in the uniquely represented communities. Notice that these constitute $\frac{1}{4}$ of all communities, so recovering a constant fraction of these vertices immediately translates to recovering a constant fraction of all vertices in the graph as well.

All the variations of our (slightly different) algorithms will be centered on Step 2 --- the computation of a weighted non-backtracking walk statistic. Thus, the final task of this subsection is to verify that these weighted non-backtracking walk statistics can be computed in polynomial time. The following proposition from~\cite{MNS:18} provides the efficiency of this step.
\begin{proposition}[{\cite[Proposition 3.3]{MNS:18}}]
    Let $N$ be the matrix $(N_{u,v})_{u,v \in V}$ encoding non-backtracking walk statistics between all pairs of vertices (whose formal definition we postpone to before Proposition \ref{prop: path count moments}). For a graph on $n$ vertices with $m$ edges and for any vector $z$, the vector $Nz$ can be computed in time $O((m+n)k)$. 
\end{proposition}

\subsection{Path counting}\label{sec: PC}
In this subsection, we recall the combinatorial path bounds from \cite[Section 4]{MNS:18} that control the number of paths of fixed length in the complete graph with conditions on the number of self-intersections. These will be a crucial ingredient in controlling the moments of the weighted path statistics in the next section. We define a path to be a sequence of vertices in the complete graph where consecutive vertices are distinct. A path is called self-avoiding if all of its vertices are distinct.  
\begin{definition}\label{def: edge types}
    Given a path $\gamma = (v_0, \ldots, v_k)$ we say that an edge $(v_i, v_{i+1})$:
    \begin{itemize}
        \item is \textit{new} if for all $j \leq i$, $v_j \neq v_{i+1}$,
        \item is \textit{old} if there is some $j < i$ such that $\{v_j, v_{j+1}\} = \{v_i, v_{i+1}\}$,
        \item is \textit{returning} otherwise. 
    \end{itemize}
    Let $k_n(\gamma), k_o(\gamma), k_r(\gamma)$ be the number of new, old, and returning edges in $\gamma$ respectively. 
\end{definition}

\begin{definition}\label{def: tangles}
    Define an $\ell$-tangle to be a graph of diameter at most $2\ell$ that contains two cycles. We say that a path $\gamma$ is \textit{$\ell$-tangle-free} if the graph $(V(\gamma), E(\gamma))$ is $\ell$-tangle-free. A path has $t$ $\ell$-tangles if $t$ is the minimal number of edges that need to be deleted to make $\gamma$ $\ell$-tangle-free. 
\end{definition}

\begin{lemma}\label{lemma: PC1}
    For any constant $C$, if $k_r \geq 1$ and $n$ is sufficiently large, then there are at most 
    \[ n^{k_n + \frac{k_r}{2} + C\log(2ek_r)} \]
    paths $\gamma$ of length at most $C \log n$, with a fixed starting and ending point, with $k_n(\gamma) = k_n$ and $k_r(\gamma) = k_r$. 
\end{lemma}

\begin{lemma}\label{lemma: PC2}
    If $k_r \geq 1$, then there are at most 
    \[ k^{5k_r + \frac{4k_rk}{\ell} + 8k_rt}n^{k_n-1} \]
    paths $\gamma$ with $t$ $\ell$-tangles that have a fixed starting and ending point, and that satisfy $k_n(\gamma) = k_n$ and $k_r(\gamma) = k_r$. 
\end{lemma}

In the second moment, we also need to control the number of pairs of paths with certain parameters, which is the content of the final lemma in this section. 

\begin{definition}
    Let $\gamma_1, \gamma_2$ be two self-avoiding paths of length $k$. An edge $(u,v)$ of $\gamma_2$ is \textit{new with respect to $\gamma_1$} if $v \not \in V(\gamma_1)$, \textit{old with respect to $\gamma_1$} if $\{u,v\}$ appears in $\gamma_1$, and \textit{returning with respect to $\gamma_1$} otherwise. Denote these quantities by $k_{n, \gamma_1}(\gamma_2), k_{o,\gamma_1}(\gamma_2)$ and $k_{r,\gamma_1}(\gamma_2)$. 
\end{definition}

\begin{lemma}\label{lemma: PC3}
    Fix vertices $u, u', v, v'$ (not necessarily distinct). There are at most 
    \[ 2(k+1)\binom{k}{k_{r,\gamma_1}}\binom{k}{k_{r,\gamma_1}+1}(2k)^{k_{r,\gamma_1}}n^{k + k_{n,\gamma_1}-1-\mathbbm{1}\{v' \not\in \{u,v\}\}} \]
    pairs $(\gamma_1, \gamma_2)$ of length-$k$ self-avoiding paths where $\gamma_1$ goes from $u$ to $v$, $\gamma_2$ goes from $u'$ to $v'$, and where $k_{n, \gamma_1}(\gamma_2) = k_{n, \gamma_1}$ and $k_{r,\gamma_1}(\gamma_2) = k_{r, \gamma_1}$. 
\end{lemma}

\begin{remark}\label{rmk: PC}
    Each of the above lemmas can be refined in the following manner, if we impose the additional restriction that exactly $r$ of the vertices in the path lie in some fixed set $U$. By the symmetry of the complete graph, each bound is multiplied by a factor of at most $\binom{k}{r}\left(\frac{|U|}{n}\right)^r \leq \left( \frac{k|U|}{n}\right)^r$ with this additional constraint.  
\end{remark}

\subsection{Weighted path statistics}\label{sec: WPS}
In this subsection we define and analyze the weighted non-backtracking walk counts between pairs of vertices. This will serve as the quantity on which we base our estimators for the community labels. 
\begin{definition}
    The following define the weights of paths:
    \begin{itemize}
        \item The weight of an edge is $W_e = \mathbbm{1}\{e \in G\} - \frac{d}{n}$.
        \item For a path $\gamma$, its weight is $X_\gamma = \prod_{e \in \gamma} W_e$. 
    \end{itemize}
\end{definition}
\begin{definition}
    Before proceeding we introduce two useful quantities.
    \begin{itemize}
        \item The first can be interpreted as the strength of the signal: $s = s_n = d\lambda = \frac{a-b}{q}$. The aforementioned KS threshold can be equivalently characterized by $s^2 = d$. 
        \item We will also make heavy use of the alignment weight of two vertices, defined as 
        \[
        a(u, v) = q\mathbbm{1}\{\sigma^\star(u) = \sigma^\star(v)\} - 1.
        \]
    \end{itemize}
\end{definition} 
The first lemma describes the moments for the weight of a self-avoiding walk or simple cycle. This will serve as the building block for our more involved computations. 
\begin{lemma}\label{lemma: weight moments}
Let $\gamma$ be a self-avoiding walk of length $k$ from $u$ to $v$. Suppose that $dmk \leq n^{1-\epsilon}$. Then
    \[ \E{X_\gamma^m \middle\vert \sigma^\star_u, \sigma^\star_v} = \begin{cases}
        \frac{a(u,v) s^k}{n^k} & m = 1, \\
        (1 + o(1))\frac{a(u,v) s^k + d^k}{n^k} & m \geq 2.
    \end{cases} \]
\end{lemma}
\begin{proof}
    We prove this by induction on $k$. In the case $k=1$, the walk is simply the edge $uv$, and it is easy to compute that 
    \[ \E{W_{uv} \middle\vert \sigma^\star_u, \sigma^\star_v} = \frac{a(u,v)s}{n}. \]
    Now suppose $\gamma$ is a self-avoiding walk of length $k$ from $u$ to $v$, and let $v'$ be the unique neighbor of $v$ in $\gamma$. We decompose based on the label of $v'$:
    \[ \E{X_\gamma \middle\vert \sigma^\star_u, \sigma^\star_v} = \E{\E{X_\gamma \middle\vert \sigma^\star_u, \sigma^\star_v, \sigma^\star_{v'}}} = \E{\E{X_{\gamma \setminus \{v'v\}}\middle\vert \sigma^\star_u,\sigma^\star_{v'}}\E{W_{v'v}\middle\vert \sigma^\star_{v'},\sigma^\star_v}}. \]
    Suppose first that $\sigma^\star_u = \sigma^\star_v$. Then the expectation over $\sigma^\star_{v'}$ factors into $\sigma^\star_{v'} = \sigma^\star_v$ and $\sigma^\star_{v'} \neq \sigma^\star_v$. By the induction hypothesis, since $\gamma \setminus \{v'v\}$ is a self-avoiding walk of length $k-1$, we obtain
    \begin{align*}
        \E{X_\gamma \middle\vert \sigma^\star_u= \sigma^\star_v} &= \frac{1}{q}\cdot \frac{(q-1)s^{k-1}}{n^{k-1}} \cdot \frac{(q-1)s}{n} + \frac{q-1}{q} \cdot \frac{-s^{k-1}}{n^{k-1}} \cdot \frac{-s}{n} = \frac{(q-1)s^k}{n^k}.
    \end{align*}
    On the other hand, now suppose $\sigma^\star_u \neq \sigma^\star_v$. The expectation over $\sigma^\star_{v'}$ now factors into three cases, where $\sigma^\star_{v'} = \sigma^\star_u$, $\sigma^\star_{v'} = \sigma^\star_v$, or $\sigma^\star_{v'}$ is different from both $\sigma^\star_u$ and $\sigma^\star_v$. Once again, inductively we have 
    \begin{align*}
        \E{X_\gamma \middle\vert \sigma^\star_u\neq \sigma^\star_v} &= \frac{1}{q}\cdot \frac{(q-1)s^{k-1}}{n^{k-1}}\cdot \frac{-s}{n} + \frac{1}{q}\cdot \frac{-s^{k-1}}{n^{k-1}} \cdot \frac{(q-1)s}{n} + \frac{q-2}{q}\cdot \frac{-s^{k-1}}{n^{k-1}} \cdot \frac{-s}{n} = \frac{-s^k}{n^k}.
    \end{align*}
    For the case $m\geq 2$, we have 
    \begin{align*}
        \E{W_{uv}^m \vert \sigma^\star_u,\sigma^\star_v} &= \left(-\frac{d}{n}\right)^m\left(1 - \frac{a(u,v)s + d}{n}\right) + \left( 1 - \frac{d}{n}\right)^m \cdot \frac{a(u,v)s+d}{n} \\
        &= \left(1 - \frac{dm}{n} + \binom{m}{2}\frac{d^2}{n^2} + \cdots + \frac{d^{m-1}m}{n^{m-1}}\right)\cdot \frac{a(u,v)s + d}{n} + \left(-\frac{d}{n}\right)^m \\
        &= \left(1 + O\big(\frac{dm}{n} + (\frac{d}{n})^m \cdot \frac{n}{b}\big)\right) \cdot \frac{a(u,v)s+d}{n}.
    \end{align*}
    Conditioning on a full labeling $\tau$ of $\gamma$ that agrees with $\sigma^\star_u$ and $\sigma^\star_v$, the above estimate gives 
    \[ \E{X_\gamma^m \middle\vert \tau} = \left(1 + O\big(\frac{dm}{n} + (\frac{d}{n})^{m-1} \cdot \frac{1}{1-\lambda}\big)\right)^{k}\prod_{xy \in \gamma} \frac{a(x,y)s+d}{n}. \]
    Notice that in the average over $\tau$, the only terms in the expansion of the product that do not vanish are $\prod_{xy \in \gamma} \frac{a(x,y)s}{n}$ and $\prod_{xy \in \gamma} \frac{d}{n}$. The first contributes exactly the first moment, and the second is identically $\frac{d^k}{n^k}$, which is the desired result. 
\end{proof}

The previous lemma allows us to estimate the weights of self-avoiding walks. In order to extend this to more general walks, we systematically decompose walks into self-avoiding pieces, and apply the lemma on each piece. We start with the definition of a valid decomposition.
\begin{definition}
    Consider a path $\gamma$. We say that a collection of paths $\zeta^{(1)}, \ldots, \zeta^{(r)}$ is a SAW-decomposition of $\gamma$ if 
    \begin{itemize}
        \item each $\zeta^{(i)}$ is a self-avoiding path,
        \item the interior vertices of each $\zeta^{(i)}$ are not contained in any other $\zeta^{(j)}$, nor is any interior vertex of $\zeta^{(i)}$ equal to the starting or ending vertex of $\gamma$,
        \item the $\zeta^{(i)}$ cover $\gamma$, meaning each edge of $\gamma$ appears at least once in some $\zeta^{(i)}$. 
    \end{itemize}
    Given a SAW-decomposition as above we let $V_\mathrm{end}$ denote the set of vertices that are endpoints of some $\zeta^{(i)}$ and we let $m_i$ denote the number of times that $\zeta^{(i)}$ was traversed in $\gamma$. 
\end{definition}

The next lemma proven in \cite{MNS:18} shows that such a decomposition exists for any walk, and describes the set of endpoints of the self-avoiding segments. We think of the set $U$ as the set of vertices for which we know the community label. Since the labels in this set are fixed, we need to ensure that the vertices appear only at the endpoints of the walks, and not in the interior.
\begin{lemma}[{\cite[Section 5.2]{MNS:18}}]\label{lemma: SAW decomp}
    For any set of vertices $U$ and any path $\gamma$ from $u$ to $v$, there exists a decomposition, which we refer to as the \textit{$U$-canonical SAW-decomposition of $\gamma$} with 
    \[ V_{\mathrm{end}} = U \cup V_{\geq 3} \cup \{u, v\} \cup \{w \in \gamma: \gamma \text{ backtracks at $w$}\} \]
    where $V_{\geq 3}$ is the set of vertices in $\gamma$ that have degree at least 3 in $\gamma$.
\end{lemma}

We now compute the expected weight of a SAW-decomposition using the building block \\ Lemma~\ref{lemma: weight moments}. These estimates on the weight of a SAW decomposition will in turn be the building blocks for our second moment estimates, which are the main goal of this subsection.
\begin{lemma}\label{lemma: decomp weight}
    Let $\zeta^{(1)}, \ldots, \zeta^{(r)}$ be a SAW-decomposition of $\gamma$. Let $z_i$ be the length of $\zeta^{(i)}$ and $e(\gamma)$ be the number of edges in $\gamma$. Then
    \[ \E{\prod_i \prod_{e \in \zeta^{(i)}} W_e^{m_i} \middle\vert \sigma^\star_{V_{\mathrm{end}}}} = (1 + o(1))\prod_{m_i = 1} \frac{a(u_i, v_i)s^{z_i}}{n^{z_i}} \prod_{m_i \geq 2} \frac{a(u_i, v_i) s^{z_i} + d^{z_i}}{n^{z_i}}. \]
    Moreover, if $r_+ = |\{i: \sigma^\star_{u_i} = \sigma^\star_{v_i}\}|$ and $e_- = \sum_{i: m_i \geq 2, \sigma^\star_{u_i} \neq \sigma^\star_{v_i}} z_i$, then 
    \[ \E{\prod_i \prod_{e \in \zeta^{(i)}} W_e^{m_i} \middle\vert \sigma^\star_{V_{\mathrm{end}}}} \leq (1+o(1))q^{r_+}s^kn^{-e(\gamma)}\left(\frac{d}{s^2}\right)^{e_-}. \]
\end{lemma}
\begin{proof}
    The first equality follows from Lemma~\ref{lemma: weight moments} and the fact that disjoint self-avoiding paths are independent conditioned on their endpoints. For the second bound, 
    \begin{align*}
        \prod_{m_i = 1} \frac{a(u_i, v_i)s^{z_i}}{n^{z_i}} \prod_{m_i \geq 2} \frac{a(u_i, v_i) s^{z_i} + d^{z_i}}{n^{z_i}} &\leq q^{r_+} \prod_{m_i = 1} \frac{s^{z_i}}{n^{z_i}} \prod_{m_i \geq 2, \sigma^\star_{u_i} = \sigma^\star_{v_i}} \frac{s^{z_i}}{n^{z_i}} \prod_{m_i \geq 2, \sigma^\star_{u_i} \neq \sigma^\star_{v_i}} \frac{d^{z_i}}{n^{z_i}} \\
        &\leq q^{r_+} s^k n^{-e(\gamma)} \prod_{e \in \zeta^{(i)}: ~m_i \geq 2, \sigma^\star_{u_i} \neq \sigma^\star_{v_i}} \frac{d}{s^{m_e}} \\
        &\leq q^{r_+} s^k n^{-e(\gamma)}\left(\frac{d}{s^2}\right)^{e_-}.
    \end{align*}
    The final inequality follows since $\frac{d}{s^{m_i}} \leq \frac{d}{s^2}$.
\end{proof}

To demonstrate the utility of the above lemmas, we begin with a simple consequence that says paths with too many returning edges contribute negligibly. The second moment estimates will follow from a similar, but much more involved computation. This estimate will also allow us to truncate the number of returning edges in our walks at certain points in our second moment computation.
\begin{lemma}\label{lemma: many intersections}
    Suppose $\frac{q|U|}{n} = o(1)$ and $k = \beta \log n$. Then for some $k^*(\chi, d, s, \beta)$, we have 
    \[ \sum_{\gamma: k_r(\gamma) > k^*} \E{X_\gamma \middle\vert \sigma^\star_U} \leq n^{-4} \]
    uniformly over assignments on $U$. 
\end{lemma}
\begin{proof}
    From Lemma \ref{lemma: decomp weight} we obtain a bound of the form 
    \[ \E{X_\gamma \vert \sigma^\star_{V_\mathrm{end}}} \leq q^{|\{i: \sigma^\star_{u_i} = \sigma^\star_{v_i}\}|} s^k n^{-e(\gamma)}\left(\frac{d}{s^2}\right)^{e(\gamma)}. \]
    Averaging over the labels of $V_\mathrm{end} \setminus U$ we obtain the bound 
    \[ \E{X_\gamma \vert \sigma^\star_U} \leq q^{|\gamma \cap U|} s^k n^{-e(\gamma)}\left(\frac{d}{s^2}\right)^{e(\gamma)}. \]
    Summing over paths with $|\gamma \cap U| = m$ using Lemma \ref{lemma: PC1}, the total weight is 
    \[ \sum_{\gamma: |\gamma \cap U| = m} \E{X_\gamma \vert \sigma^\star_U} \leq \frac{q^{m+1}}{n^{m(1+o(1))}} s^k \left(\frac{d}{ns^2}\right)^{e(\gamma)} \cdot n^{k_n + k_r/2 + O(\log k_r)}. \]
    The sum over $m$ converges, so 
    \[ \sum_{\gamma: k_r(\gamma)>k^*} \E{X_\gamma \vert \sigma^\star_U} \leq s^k \left(\frac{d}{s^2}\right)^{k_n+k_r} \cdot n^{-k_r/2 + O(\log k_r)}. \]
    For $k_r$ sufficiently large we have a total contribution of at most $n^{-4}$. 
\end{proof}

The next two propositions are the main technical components of our work, which describe the first two moments of the weighted path counts along with various correlations. The form of the estimates depends on the regime of parameters we are considering. First, below the KS threshold, but above the modified KS threshold, we have the following moments. These will be essential in proving Theorem \ref{thm: main alg}(1). The random variables we are interested in are the following, for a fixed set of vertices $U$.
\begin{itemize}
    \item For any two fixed vertices $u,v$, the \textit{weighted self-avoiding walk count} between $u$ and $v$ is $S_{u,v}^U = \sum_{\gamma} X_\gamma$ where $\gamma$ sums over all self-avoiding walks of length $k$ from $u$ to $v$ that avoid $U$, aside from possibly $u$ and $v$.
    \item The \textit{weighted non-backtracking walk count} between $u$ and $v$ is $N_{u,v}^U = \sum_{\gamma} X_\gamma$ where $\gamma$ now sums over all non-backtracking walks of length $k$ from $u$ to $v$, again avoiding $U$ aside from $u$ or $v$.
\end{itemize}

In the next two propositions we also take a slightly more general setup than the block model. As alluded to in Subsection~\ref{sec: Alg} our algorithm requires some preprocessing of our sample from $\SBM(n, q, d, \lambda)$ involving deleting some vertices and choosing representatives. This slightly biases the community distribution of the vertices in the remaining part of the graph. However, and importantly, conditioned on the communities edges are still drawn independently and at random with probabilities $\frac{a}{n}$ and $\frac{b}{n}$. Thus, we take a setup in which community sizes are fixed as $(1+o(1))\frac{n}{q}$ and edges are drawn according to the $\SBM$ probabilities conditioned on these labels. 
\begin{definition}
    Suppose that $n - n' = o(\frac{n}{\log n})$, $\abs{s_i - \frac{n}{q}} = o(\frac{n}{q\log n})$ and $\sum_{i=1}^q s_i = n$. Let $G \sim \widetilde{\SBM}(n', (s_i)_{i=1}^q, \frac{a}{n}, \frac{b}{n})$ as follows. Partition $V(G)$ uniformly at random into communities $C_1, \ldots, C_q$ such that $|V_i| = s_i$. For any two vertices $u, v$ independently include the edge $uv$ with probability $\frac{a}{n}$ if $u$ and $v$ are in the same community, and probability $\frac{b}{n}$ otherwise.
\end{definition}

\begin{proposition}\label{prop: path count moments}
    Let $G \sim \widetilde{\SBM}(n', (s_i)_{i=1}^q, \frac{a}{n}, \frac{b}{n})$. Suppose $q = o(\frac{n}{\log^3 n})$ and that $d\lambda^2 < 1$ but $d\lambda^{1/\chi} > 1$. Choose $\beta$ such that $\frac{1-\chi}{\log s} < \beta < \frac{2\chi-1}{\log (d/s^2)}$. Let $u, v, u', v'$ be distinct vertices. Set $k = \lfloor \beta \log n \rfloor$. Suppose that $U, U' \subseteq V$ contain $\{u,v\}$ and $\{u', v'\}$ respectively, and that both have cardinality at most $o(\frac{n}{\log n})$. Then the following hold:
    \begin{align}
        \label{eq: first moment}&\E{S_{u,v}^U \middle\vert \sigma^\star_U} = (1 + o(1)) \frac{a(u,v) s^k}{n}, \\
        \label{eq: second moment}&\E{(S_{u,v}^U)^2 \middle\vert \sigma^\star_U} \leq \begin{cases}
            (1 + o(1)) \frac{s^2}{(s-1)^2} \cdot \frac{q^2 s^{2k}}{n^2} & \sigma^\star_u = \sigma^\star_v, \\
            (1+o(1))\frac{d^k}{n} & \sigma^\star_u \neq \sigma^\star_v,
        \end{cases} \\
        \label{eq: cross moment}&\E{S_{u,v}^{U \cup V}S_{u',v'}^{U \cup U'} \middle\vert \sigma^\star_U, \sigma^\star_{U'}} = (1 + O(q/n))\E{S_{u,v}^{U \cup U'}\middle\vert \sigma^\star_U, \sigma^\star_{U'}}\E{S_{u',v'}^{U \cup U'}\middle\vert \sigma^\star_U, \sigma^\star_{U'}}, \\
        \label{eq: complicated paths}&\Prob{\abs{S_{u,v}^U - N_{u,v}^U} > \epsilon qs^kn^{-1}} \leq \begin{cases}
            \epsilon^{-2}qn^{-1} & \sigma^\star_u = \sigma^\star_v, \\
            \epsilon^{-2}n^{-1} & \sigma^\star_u \neq \sigma^\star_v.
        \end{cases}
    \end{align}
\end{proposition}
The next proposition describes the analogous moments above the KS threshold. The key difference in this setting is the estimate on the second moment. These will be essential in proving Theorem \ref{thm: main alg}(2).
\begin{proposition}\label{prop: path count moments KS}
    Let $G \sim \widetilde{\SBM}(n', (s_i)_{i=1}^q, \frac{a}{n}, \frac{b}{n})$. Suppose that $q = o(\frac{n}{\log^3 n})$ and  $d\lambda^2 > 1$. Choose $\beta$ such that $\frac{2}{\log (s^2/d)} < \beta$. Let $u,v,u',v'$ be distinct vertices. Set $k = \lfloor \beta \log n \rfloor$. Suppose that $U, U' \subseteq V$ contain $\{u,v\}$ and $\{u',v'\}$ respectively, and that both have cardinality at most $o(\frac{n}{\log n})$. Then (\ref{eq: first moment}), (\ref{eq: cross moment}), (\ref{eq: complicated paths}) and the following hold:
        \begin{equation}\label{eq: second moment 2}
            \E{(S_{u,v}^U)^2 \middle\vert \sigma^\star_U} \leq \begin{cases}
            (1 + o(1)) \frac{s^2}{(s-1)^2} \cdot \frac{q^2 s^{2k}}{n^2} & \sigma^\star_u = \sigma^\star_v, \\
            (1+o(1)) \frac{s^2}{s^2-d} \cdot \frac{qs^{2k}}{n^2} & \sigma^\star_u \neq \sigma^\star_v.
            \end{cases}
        \end{equation}
\end{proposition}
The proof of these two propositions is largely the same in terms of bookkeeping the powers of $q$. For this reason, we provide the two proofs in parallel and preview the  differences here. The key difference arises when terms of the form $d^k$ are all asymptotically dominated by terms of the form $s^{2k}$ above the KS threshold, whereas below the KS threshold the converse is true. In particular, the terms $\left(\frac{d}{s^2}\right)^{e_-}$ in the conclusion of Lemma~\ref{lemma: decomp weight} contribute a factor smaller than 1 in the setting of Proposition~\ref{prop: path count moments KS}.

We emphasize here that when $G \sim \widetilde{\SBM}(n', (s_i)_{i=1}^q, \frac{a}{n}, \frac{b}{n})$ compared to $\SBM(n,q,a,b)$, the community labels are drawn without replacement from a slightly biased distribution rather than i.i.d.\ uniformly at random. However, the bias of the distribution of sampling $k$ labels can be quantified as follows:
\begin{align*}
    \prod_{j=1}^k \frac{s_{i_j}-j+1}{n'-j+1} &= \prod_{j=1}^k \left(\frac{s_{i_j}}{n'-j+1} -\frac{j-1}{n'-j+1}\right) \prod_{j=1}^k \left(1 + o\left(\frac{1}{\log n}\right)\right) \cdot \left(\frac{1}{q} + \frac{\log n}{n} \right) \\
    &= \left( 1 + o\left(\frac{1}{\log n}\right)\right)^{k} \cdot \frac{1}{q^k} = (1+o(1))\frac{1}{q^k}.
\end{align*}
for $k = O(\log n)$. Since we consider only walks of length $O(\log n)$ the effect is negligible. Similarly, since the biased degree $d' = (1 + o(\frac{1}{\log n}))d$ and $n' = (1 + o(\frac{1}{\log n}))n$ the quantities $(d')^k$ and $(n')^k$ are also $(1+o(1))d^k$ and $(1+o(1))n^k$ respectively for $k = O(\log n)$. 

For example, in the proof of Lemma~\ref{lemma: weight moments} replacing $\frac{1}{q}$ by $\frac{1+o(1)}{q}$ costs us only a $1+o(1)$ in the estimates. Thus, in the following proofs, we will refer only to $d$ and $n$ with the understanding that the $1+o(1)$ factors are absorbed into the leading constants. In essence, the moments are computed according to $\SBM(n,q,a,b)$, but the result holds for $\widetilde{\SBM}(n', (s_i)_{i=1}^q, \frac{a}{n}, \frac{b}{n})$ as the bias is negligible.

\begin{proof}[Proof of (\ref{eq: first moment})]
    The number of walks which avoid $U$ is 
    \begin{multline*}
        (n-|U|)(n-|U|-1) \cdots (n-|U|-k+1) = n^k\left(1 - \frac{|U|}{n}\right)\cdots \left(1-\frac{|U|+k-1}{n}\right) \\ \geq n^k\exp(-2\sum_{i=1}^k \frac{|U|+i-1}{n}) = n^k\exp(- \frac{2k|U|}{n} - \frac{2k^2}{n}) = (1+o(1))n^k.
    \end{multline*}
    The claim follows from Lemma~\ref{lemma: weight moments} in the settings of both Propositions~\ref{prop: path count moments} and~\ref{prop: path count moments KS}.
\end{proof}
\begin{proof}[Proof of (\ref{eq: second moment}) and (\ref{eq: second moment 2})]
    We control the paths in the 4 groups based on their number of returning edges: $\Gamma_0 = \{k_r = 0\}$, $\Gamma_1 = \{k_r = 1\}$, and $\Gamma_2 = \{2 \leq k_r \}$. For this and the following proofs, we crudely upper bound the number of paths avoiding $U$ by the total number of paths.
    
    For $\gamma \in \Gamma_0$, the total weight of paths conditioned on $\sigma^\star_U$ is 
    \[ \E{X_\gamma \middle\vert \sigma^\star_U} = (1+o(1))\frac{a(u,v)s^k + d^k}{n} \leq \begin{cases}
        (1+o(1))q \cdot \frac{s^k}{n} & \sigma^\star_u = \sigma^\star_v, \\
        (1+o(1)) \cdot \frac{d^k}{n} & \sigma^\star_u \neq \sigma^\star_v.
    \end{cases} \]
    By choice of $\beta$, we have $n \ll qs^k$, which implies that $\E{X_\gamma \middle\vert \sigma^\star_U} \ll \left( \frac{qs^k}{n}\right)^2$ when $\sigma^\star_u = \sigma^\star_v$. This contributes the main term to (\ref{eq: second moment}) when $\sigma^\star_u \neq \sigma^\star_v$, but is dominated by $\frac{qs^{2k}}{n^2}$ in the setting of (\ref{eq: second moment 2}).

    For $\gamma \in \Gamma_1$, the paths $\gamma$ consist of a cycle with two tails. Let $u'$ and $v'$ denote the vertices where the tails containing $u$ and $v$ respectively attach to the cycle. We decompose according to $i$, the distance from $u$ to $u'$, and $j$, the distance from $v$ to $v'$. Then we have 
    \[ \E{X_\gamma \middle\vert \sigma^\star_U} = \sum_{i,j \leq k}(1+o(1))\left(a(u,u')s^{i}+d^{i}\right) \cdot \left(a(v,v')s^{j}+d^{j}\right) \cdot a(u',v')^2s^{2k-2(i+j)} \cdot n^{-2}. \]
    If $\sigma^\star_u = \sigma^\star_v$, the main term is contributed by $\sigma^\star_u = \sigma^\star_{u'} = \sigma^\star_{v'} = \sigma^\star_v$ where the $a(\cdot, \cdot)$ terms contribute $q^4$, which occurs with $\frac{1 + o(1)}{q^2}$ probability. When $\sigma^\star_{u'} = \sigma^\star_{v'} \neq \sigma^\star_u$, the $a(\cdot ,\cdot)$ terms contribute $q^2$, and this occurs with probability $\frac{1 + o(1)}{q}$. When $\sigma^\star_u = \sigma^\star_{u'} \neq \sigma^\star_{v'}$ the $a(\cdot, \cdot)$ terms contribute $q$, while the probability this occurs is only $\frac{1+o(1)}{q}$. Finally, in the last remaining case in which $\sigma^\star_u$, $\sigma^\star_{u'}$, and $\sigma^\star_{v'}$ are all different, the $a(\cdot ,\cdot)$ terms contribute only 1. In this case analysis, we crucially use the fact that the contribution from a single factor of $q$ dominates the contribution from $d^k$ when compared to $s^k$. Thus, all together in this situation we have the bound 
    \[ \E{X_\gamma \middle\vert \sigma^\star_U} \leq \frac{1 + o(1)}{q^2} \sum_{i,j \leq k} \frac{q^4s^{2k-i-j}}{n^2} \leq (1 + o(1))\frac{s^2}{(s-1)^2} \cdot \left( \frac{qs^k}{n}\right)^2. \]
    
    If $\sigma^\star_u \neq \sigma^\star_v$, a similar case analysis yields that the main term is contributed by the case $\sigma^\star_{u'} = \sigma^\star_{v'}$, which occurs with probability $\frac{1+o(1)}{q}$. In the setting of (\ref{eq: second moment}) we have the bound 
    \[ \E{X_\gamma \middle\vert \sigma^\star_U} \leq \frac{1+o(1)}{q}\sum_{i+j \leq k} \frac{q^2s^{2k}}{n^2}\left( \frac{d}{s^2}\right)^{i+j} \leq k \cdot \frac{\left(\frac{d}{s^2}\right)^{k+1} - 1}{\frac{d}{s^2} - 1} \cdot \frac{qs^{2k}}{n^2} \lesssim \left(\frac{d}{s^2}\right)^{k}\cdot \frac{qs^{2k}}{n^2} \ll \frac{d^k}{n}. \]
    In the setting of (\ref{eq: second moment 2}), the factor $\frac{d}{s^2} < 1$, and so the bound is instead of the form 
    \[ \E{X_\gamma \middle\vert \sigma^\star_U} \leq \frac{1+o(1)}{q}\sum_{i+j \leq k} \frac{q^2s^{2k}}{n^2}\left( \frac{d}{s^2}\right)^{i+j} \leq k \cdot \frac{1 - \left(\frac{d}{s^2}\right)^{k+1}}{1 - \frac{d}{s^2}} \cdot \frac{qs^{2k}}{n^2} \lesssim \frac{s^2}{s^2-d} \cdot \frac{qs^{2k}}{n^2}. \]
    This constitutes the main contribution. 

    For $\gamma \in \Gamma_2$, we further partition based on $k_n$. By Lemma~\ref{lemma: PC3} the number of such paths is at most $(2k)^{4k_r}n^{k_n-1}$. By Lemma~\ref{lemma: decomp weight}, conditioned on $\sigma^\star_{V_\mathrm{end}}$ the contribution of such paths is bounded by $q^{r_+}s^{2k}n^{-k_n-k_r}\left(\frac{d}{s^2}\right)^{e_-}$. In the average over labels of $V_\mathrm{end} \setminus U$ we have a few cases. 
    
    If $\sigma^\star_u \neq \sigma^\star_v$ the exponent of $q$ is at most $k_r$ and the exponent of $\frac{d}{s^2}$ is at most $k$. The inequality $e_- \leq k$ follows since $e_-$ counts a subset of edges that are traversed at least twice in the walk $\gamma$, which is of total length $2k$. Thus, there can be at most $k$ of these edges. In the average over $\sigma^\star_{V_\mathrm{end} \setminus U}$, the exponent of $q$ is at most $k_r$ by the following argument. Each returning edge that has an endpoint in $V_\mathrm{end}$ introduces at most one factor of $q$ on average, by an analogous argument to the $\Gamma_1$ case above. Since $\sigma^\star_u \neq \sigma^\star_v$, we start with no factors of $q$, and thus the total expected contribution is at most $q^{k_r}$. This yields a bound of $q^{k_r}d^k n^{-k_n - k_r}$. This is a total contribution of 
    \begin{align*}
        \sum_{k_r \geq 2} \sum_{k_n \leq k} (2k)^{4k_r}n^{k_n-1} \cdot q^{k_r + r}d^{k}n^{-k_n-k_r} &= \frac{d^{k}}{n}\sum_{k_r \geq 2} \left(\frac{q}{n}\right)^{(-1+o(1))k_r} \lesssim \frac{d^k}{n} \cdot \frac{q^2}{n^2} \ll \frac{d^{k}}{n}.
    \end{align*} 
    In the setting of (\ref{eq: second moment 2}), the factors of $\frac{d}{s^2}$ can be dropped, and the same computation yields the bound $\frac{s^{2k}}{n} \cdot \frac{q^2}{n^2} \ll \frac{qs^{2k}}{n^2}.$ 
    
    If $\sigma^\star_u = \sigma^\star_v$, by the argument in the previous case, and the fact that $\sigma^\star_u = \sigma^\star_v$ gives an additional factor of $q$ to start, the exponent of $q$ under the averaging is at most $k_r + 1$. A key observation is that the argument implies that this occurs precisely when all labels on $V_\mathrm{end}$ are the same. In particular, we can have no contributions of the form $\frac{d}{s^2}$, so the total contribution is bounded by $q^{k_r + 1}s^{2k}n^{-k_n - k_r}$. On the other hand, if the labels are not all the same, we may have up to $k$ factors of $\frac{d}{s^2}$ as argued above, but lose the extra factor of $q$. Thus, the contribution is bounded by $q^{k_r} d^k n^{-k_n - k_r}$. Since $\left(\frac{d}{s^2}\right)^k \ll q$ this is a smaller contribution. 
    All together, we have a bound of the form
    \begin{align*}
        \sum_{k_r \geq 2} \sum_{k_n \leq k} (2k)^{4k_r}n^{k_n-1}\cdot q^{k_r + r + 1}s^{2k}n^{-k_n-k_r} &= \frac{qs^{2k}}{n}\sum_{k_r \geq 2} \left(\frac{q}{n}\right)^{(-1+o(1))k_r} \lesssim q^3s^{2k}n^{-3 + o(1)} \ll \frac{q^2s^{2k}}{n^2}.
    \end{align*}    
\end{proof}
\begin{proof}[Proof of (\ref{eq: cross moment})]
    The pairs of paths are once again divided into groups based on the number of intersections: $\Gamma_0 = \{ k_{r, \gamma_1}(\gamma_2) = 0\}$, and $\Gamma_1 = \{ k_{r, \gamma_1}(\gamma_2) \geq 1 \}$. $\Gamma_0$ contributes the main term, in which the weights $X_{\gamma_0}$ and $X_{\gamma_1}$ are truly independent. 

    For $\Gamma_1$, we proceed by partitioning based on $k_{n, \gamma_1}(\gamma_2)$. By Lemma~\ref{lemma: PC3}, the number of such pairs is $n^{k + k_{n, \gamma_1} - 2 + o(1)}$. When both $\sigma^\star_u = \sigma^\star_v$ and $\sigma^\star_{u'} = \sigma^\star_{v'}$, we have that the contribution from $\Gamma_0$ is of order $\frac{q^2s^{2k}}{n^2}$. In this case, by Lemma~\ref{lemma: decomp weight} and similar reasoning to the proof of (\ref{eq: second moment}), the contributions are bounded by $q^{k_{r, \gamma_1} + 2}s^{2k}n^{-k - k_{n, \gamma_1}-k_{r, \gamma_1}}$. Here the $+2$ in the exponent of $q$ is because we start with two walks, which can contribute two factors of $q$ initially. In total this contributes
    \begin{multline*}
        \sum_{k_{r, \gamma_1} \geq 1} \sum_{k_{n, \gamma_1} \leq k} n^{k + k_{n, \gamma_1} - 2 + o(1)} \cdot q^{k_{r, \gamma_1} + r + 2}s^{2k}n^{-k - k_{n, \gamma_1}-k_{r, \gamma_1}} = \frac{q^2s^{2k}}{n^2}\sum_{k_{r, \gamma_1} \geq 1} \left(\frac{q}{n}\right)^{k_{r, \gamma_1}} \ll \frac{q^2s^{2k}}{n^2}.
    \end{multline*}
    When $\sigma^\star_u = \sigma^\star_v$ and $\sigma^\star_{u'} \neq \sigma^\star_{v'}$, we have that the contribution from $\Gamma_0$ is of order $\frac{qs^{2k}}{n^2}$. In this case we again have one initial factor of $q$ from the labels of the endpoints, so by Lemma~\ref{lemma: decomp weight} the contributions are bounded by $q^{k_{r, \gamma_1} + 1}s^{2k}n^{-k - k_{n, \gamma_1}-k_{r, \gamma_1}}$, which gives a total contribution of $\ll \frac{qs^{2k}}{n^2}$. 
    Finally when $\sigma^\star_u \neq \sigma^\star_v$ and $\sigma^\star_{u'} \neq \sigma^\star_{v'}$ the contribution from $\Gamma_0$ is of order $\frac{s^{2k}}{n^2}$. In this case the contributions from each pair of paths are bounded by $q^{k_{r, \gamma_1}}s^{2k}n^{-k - k_{n, \gamma_1}-k_{r, \gamma_1}}$, and so the total contribution is $\ll \frac{s^{2k}}{n^2}$. All together, this confirms that the independent case $\Gamma_0$ is indeed the main contribution.
\end{proof}
\begin{proof}[Proof of (\ref{eq: complicated paths})]
    With the possible addition of SAW-decomposition components which are cycles, we have that the contribution from each individual path is bounded by $q^{k_{r} + 1 + \mathbbm{1}\{\sigma^\star_u = \sigma^\star_v\}}s^{2k}n^{-k_{n}-k_{r}}$ (recall in the self-avoiding walk case we had an exponent of $k_{r} + \mathbbm{1}\{\sigma^\star_u = \sigma^\star_v\}$ instead). Note that since $\gamma$ is the concatenation of two walks that are non-backtracking but not self-avoiding, we must have $k_r(\gamma) \geq 1$. Thus, Lemma \ref{lemma: PC2} implies that the number of such paths is at most $n^{k_n(\gamma_1) + k_n(\gamma_2)-2 + o(1)}$. Thus, the total contribution is at most 
    \begin{multline*}
        \sum_{k_r \geq 1} n^{k_n -2 + o(1)} \cdot q^{k_{r} + r + 1 + \mathbbm{1}\{\sigma^\star_u = \sigma^\star_v\}}s^{2k}n^{-k_{n}-k_{r}} = \\ \frac{q^{1 + \mathbbm{1}\{\sigma^\star_u = \sigma^\star_v\}}s^{2k}}{n^2} \sum_{k_r \geq 1}\left( \frac{q}{n}\right)^{-(1 + o(1))k_r} 
        \lesssim \frac{q^{2 + \mathbbm{1}\{\sigma^\star_u = \sigma^\star_v\}}s^{2k}}{n^3}.
    \end{multline*}
    By Chebyshev's inequality, 
    \begin{align*}
        \Prob{|S_{u,v} - N_{u,v}| > \epsilon qs^kn^{-1}} &\leq \frac{q^{2 + \mathbbm{1}\{\sigma^\star_u = \sigma^\star_v\}}s^{2k}n^{-3}}{\epsilon^2q^2s^{2k}n^{-2}} = \epsilon^{-2}q^{\mathbbm{1}\{\sigma^\star_u = \sigma^\star_v\}}n^{-1}.
    \end{align*}
\end{proof}

\subsection{Below the KS threshold -- Theorem \ref{thm: main alg}(1)}\label{sec: proof thm WR} 
\begin{algorithm}[]\label{alg: below KS}
    \SetAlgoLined
    Take $U \subseteq V$ to be a random subset of size $O(n^{\frac{\chi+1}{2}})$.  \\
    Draw $q$ vertices $U_* = \{u_1, \ldots, u_q\}$ uniformly at random from all vertices in $U$ having degree at least $\log\log n$ to $V\setminus U$. \\
    \For{$w \in V \setminus U$}{
    Partition $V \setminus (U \cup \{w\})$ randomly into $M = \Theta(\log d)$ nearly equal sized parts $V_1(w), \ldots, V_M(w)$, such that $\abs{|V_i(w)| - |V_j(w)|} \leq 1$ for all $i, j$. \\
    \For{$\ell \in [q]$}{
    Compute $Z_i(w, u_\ell) = \sum_{u \in N(u_\ell) \cap V_i}\sum_{v \in S_i(w)} N^{(i)}_{u,v}$ where $S_i(w) = \{v: d_{V_i}(w,v) = L\}$ and $L$ is growing sufficiently slowly for each $1 \leq i \leq M$. \\
    }
    Set $\widehat{\sigma}_w = \ell$ if $\ell$ is the unique label such that $Z_i(w, u_\ell) > 3aq\lambda\frac{s^{k+L}}{n M^{k+L-1}}|N(u_\ell) \cap V_i| - \frac{s^k}{nM^{k-1}}|N(u_\ell) \cap V_i||S_i(w)|$ for every $i$, where $a$ is as in Lemma \ref{lemma: branching process} if it exists. Otherwise assign $\widehat{\sigma_w}$ uniformly at random. 
    }
    \For{ $w \in U \setminus U_*$ } {
    Assign $\widehat\sigma(w)$ uniformly at random.
    }
    \For{$\ell \in [q]$}{
    Assign $\widehat\sigma(u_\ell) = \ell$.
    }
    \caption{Below the KS threshold}
\end{algorithm}
Armed with the results of the previous two subsections, we can finally present Algorithm~\ref{alg: below KS} and prove its correctness. Recall from Subsection~\ref{sec: Alg} that Step~1 chooses a set of representatives $U_*$, in which at least $\frac{q}{4}$ of the communities are uniquely represented and do not appear too often in each other's neighborhoods. Let $\mathcal{C}_\mathrm{good} \subseteq [q]$ denote this subset of good communities, i.e. those that are uniquely represented and appear no more than $\log n$ times along neighbors of other representatives. We detail how to carry out Step~2 of basic outline --- namely, how to decide candidates based on the walk statistics. In this regime, we perform a random partition of $V \setminus (U \cup \{w\})$ into $M = \Theta(\log d)$ disjoint vertex parts $V_1(w), \ldots, V_M(w)$. We will often write $V_i$ rather than $V_i(w)$ since we classify only a single $w$ at a time. We then compute 
\[ Z_i(w, u_*) \vcentcolon= \sum_{u \in N(u_*)\cap V_i} \sum_{v \in S_i(w)} N_{u,v}^{V_i^c}, \]
where $S_i(w) = \{v: d_{V_i}(x, v) = L\}$ for $L$ a slowly growing function of $n$ to be specified later and we recall that $N_{u,v}^{V_i^c}$ is the weighted non-backtracking walk count between $u$ and $v$ which avoid $V_i^c$. In particular, this is a statistic using only paths contained in a fixed part $V_i$ indicating whether or not $w$ aligns with $u_*$. The criteria for assigning $u_*$ as a candidate for $w$ is that $Z_i(w, u_*)$ is large \textit{for every} $i$. For the remainder of this subsection, since we restrict to quantities $N_{u,w}^{V_i^c}$ contained in a specific part we work in a model where $n_i \mapsto \frac{n}{M}$ but $a$ and $b$ remain fixed. This affects the alternate parameterization in the following way: $d_i \mapsto \frac{d}{M}$, $\lambda \mapsto \lambda$. We also write $(i)$ in the superscript rather than $V_i^c$ for ease of notation.

The reason for this modification is that the second moment does not imply sufficient concentration for the candidate list of a vertex $w$ to be small. In particular, it is only sufficiently concentrated to produce a better than random alignment, but not the desired constant alignment. Partitioning the graph into pieces is a ``boosting'' step to manually bootstrap the statistic and obtain tighter concentration. This will result in smaller candidate lists, and a constant alignment in the final output partition.

We prove that the output of the algorithm, which we denote by $\widehat \sigma$, has the desired alignment in two steps. First, we show that a local estimate based on only  neighborhoods around vertices has constant alignment with the true communities with high probability. The local predictions are defined as follows. Define 
\[ Z_i^{\mathrm{loc}}(w, u_*) \vcentcolon= \sum_{u \in N(u_*) \cap V_i}\sum_{v \in S_i(w)} a(u,v).\]
where $S_i(w) = \{v: d_{V_i}(w,v) = L\}$.
Let $\tau_w = u_*$ if 
\[ Z_i^\mathrm{loc}(w, u_*) > 4aq\lambda s_i^L|N(u_*) \cap V_i| - |N(u_*) \cap V_i||S_i(w)| \qquad \forall \ 1 \leq i \leq M \] 
and 
\[ Z_i^\mathrm{loc}(w, u) < 2aq\lambda s_i^L|N(u_*) \cap V_i| - |N(u_*) \cap V_i||S_i(w)| \qquad \forall \ 1 \leq i \leq M, \, u \neq u_* \in U_*.\] 
Otherwise assign it randomly using the same randomness that was used to generate random $\widehat{\sigma}_w$ in Algorithm~\ref{alg: below KS}. We first claim that Proposition~\ref{prop: path count moments} applies to our setup.

\begin{claim}\label{claim}
    Fix $w \in V\setminus U$. With high probability $G\vert_{V_i(w)}$ satisfies the conditions of \\ $\widetilde{\SBM}(n', (s_j)_{j=1}^q, \frac{a}{n_i}, \frac{b}{n_i})$. Moreover, given $\sigma_{U \cup w}$ the graphs $G|_{V_i(w)}$ are conditionally independent for distinct $i$.
\end{claim}
\begin{proof}
    Sample the block model $G$ as follows. First reveal the \textit{sizes} of each community. By concentration of binomial random variables, with high probability these will be $\frac{n}{q} + O(\sqrt{\frac{n}{q}})$. Then sample $U$ of size $O(n^{\frac{\chi+1}{2}}) = O(\sqrt{nq})$ and reveal $\sigma^\star_U$ and the edges between $U$ and $V \setminus U$. Once again by concentration, with high probability the community sizes on $U$ will be $\frac{|U|}{q} + O(\sqrt{\frac{|U|}{q}})$. Choose $U_*$ to be the $q$ vertices of highest degree. Finally reveal $\sigma_{N(U_*)}$ from the conditional distribution. 
    
    By a calculation similar to Lemma~\ref{lemma: init}, with high probability there will be at most $\log n\log\log n$ vertices in $U_*$ with the same label, and thus at most $\log^2 n$ vertices with the same label in $N(U_*)$. In particular, with high probability the community sizes on $V \setminus (U \cup N(U_*))$ determined by the revealed information will all be $\frac{n-|U|}{q} + O(\sqrt{\frac{n}{q}} + \sqrt{\frac{|U|}{q}} + \log^2 n)$. We conclude that 
    \[ n' = n - |U| = n - o\left(\frac{n}{\log n}\right) \text{ and } \abs{s_j - \frac{n}{q}} = O\left(\sqrt{\frac{n}{q}} + \log^2 n\right) = o\left(\frac{n}{q\log n}\right). \]

    Conditional independence follows from the independence in the original graph and the fact that $U \cup \{w\}$ separates the pieces of the partition (as we ignore all edges between distinct pieces). 
\end{proof}

With this in hand, we proceed to analyze the estimator based on the non-backtracking walks. To do this we introduce a branching process that describes the behavior on the neighborhood around a random vertex $w$. Let $\mathrm{GW}_L$ be a labeled tree $(T, \tau)$ distributed as follows. Start with a root vertex $\rho$ with label $\tau(\rho)$ uniformly random in $[q]$. Let $\rho$ have $C_\rho$ children where $C_\rho \sim \mathrm{Bin}(n, \frac{d}{n})$. For each child $c$, let $\tau(c) = \tau(\rho)$ with probability $\lambda$ and uniformly random otherwise. Iterate this procedure for depth $L$. We will take $L$ to be a sufficiently slowly growing function of $n$, e.g. $L = \lceil \log^{(10)} n \rceil$ the 10-fold iterated logarithm of $n$ suffices.

This labeled tree can be equivalently generated as follows: generate an unlabeled tree with offspring distribution $\mathrm{Bin}(n, \frac{d}{n})$ to depth $L$. Perform bond percolation on the tree by retaining each edge independently with probability $\lambda$. For each connected component in the percolated tree, label every vertex in the component with the same uniformly random label. We use this formulation for the following lemma.
\begin{lemma}\label{lemma: branching process}
    Suppose $d\lambda > 1$ and generate $T \sim \mathrm{GW}_L$. Let $N_L$ be the number of vertices at depth $L$ in the connected component of $\rho$, and let $M_L$ be the maximum number of vertices at depth $L$ of a cluster not containing $\rho$. Then there exist constants $a, \delta > 0$ depending on $d, \lambda$ but not $L$ such that
    \[ \Prob{N_L > 8a(d\lambda)^L, M_L < a(d\lambda)^L} > \delta. \]
\end{lemma}
\begin{proof}
Throughout the proof we rely on the classical fact that the size of a branching process at depth $L$ normalized by its mean is a martingale and converges almost surely to a non-degenerate limiting random variable. We condition on the positive probability event that the root has $D$ children (to be chosen later in terms of $d$ and $\lambda$) and none of them are disconnected from the root in the percolation process. Then
\[ N_L = \sum_{i=1}^{D} N_{L-1}^{(i)} \qquad \text{ and } \qquad M_L = \max_{1 \leq i \leq D} M_{L-1}^{(i)} \]
where $N_{L-1}^{(i)}$ and $M_{L-1}^{(i)}$ are i.i.d. copies of the random variables for depth $L-1$ processes. Since $\E{N_{L-1}} = (d\lambda)^{L-1}$ and $\E{N_{L-1}^2} = C(d)(d\lambda)^{2(L-1)}$, the Paley--Zygmund inequality implies 
\[ \Prob{N_L > \frac{1}{2}D(d\lambda)^{L-1}} \geq \frac{1}{4C(d)}. \]
On the other hand,
\[ \Prob{M_L > \frac{1}{16} D(d\lambda)^{L-1}} \leq  D\cdot \Prob{M_{L-1} > \frac{1}{16}D(d\lambda)^{L-1}}. \]
We bound the probability that $M_{L-1}$ is large as follows. We associate each cluster $C$ with its root $\rho_C$, that is the vertex of minimum depth. There are in expectation $d^h$ candidates for $\rho_C$ at depth $h$, and each one contributes a cluster with size independently distributed as $N_{L-h}$. Thus, for any finite moment $k$,
\[ \Prob{M_{L-1} > t} \leq \sum_{h \geq 1} d^h\Prob{N_{L-h} > t} \leq \sum_{h\geq 1} d^h \cdot \E{N_{L-h}^k} t^{-k} = \sum_{h\geq 1} d^h \cdot C_1(d)(d\lambda)^{k(L-h)}t^{-k}. \]
Setting $t = \frac{1}{16}D(d\lambda)^{L-1}$ we obtain
\[ \Prob{M_{L-1} > t} \leq \frac{16^kC_1(d)}{D^k}\sum_{h \geq 1} d^h(d\lambda)^{-k(h-1)}. \]
For $k = 2\log_{d\lambda} d$ the series on the right hand side converges. Altogether we obtain 
\[ \Prob{M_L > \frac{1}{16}D(d\lambda)^{L-1}} \leq C_2(d, \lambda)D^{-(k-1)}. \]
Thus, for $D$ sufficiently large in terms of $d$ and $\lambda$, we have the desired conclusion with $a = \frac{D}{16d\lambda}$ and $\delta = \frac{1}{8C(d)}$. 
\end{proof}

\begin{lemma}\label{lemma: coupling}
    Let $m = \lfloor \log n \rfloor$ and $v_1, \ldots, v_m$ be uniformly randomly selected (without replacement) vertices from $V$. Let $(T_i, \tau_i)_{i=1}^{m}$ be i.i.d.\ labeled trees distributed as $\mathrm{GW}_L.$ Then
    \[ d_{\mathrm{TV}}((S(v_i), \sigma^\star_{S(v_i)})_{i=1}^{m}, (T_i, \tau_i)_{i = 1}^{m}) \to 0. \]
\end{lemma}
\begin{proof}
    The proof is a small variant of a classic result in the theory of sparse random graphs. Note that with high probability, the sizes of the neighborhoods $S(v_i)$ are all $n^{o(1)}$. Thus, the probability that they are not mutually disjoint is $n^{-1+o(1)}$. Conditioned on the event that all neighborhoods are disjoint, we couple the two processes by simultaneously revealing the number of neighbors of successive vertices in $S(v_i)$ and the number of children of the corresponding vertex in $T_i$. The number of neighbors of a vertex in $S(v_i)$ is distributed as $\mathrm{Bin}(n-n^{o(1)}, \frac{d}{n})$ and the number of children of a vertex in $T_i$ is $\mathrm{Bin}(n, \frac{d}{n})$. Since $d_{\mathrm{TV}}(\mathrm{Bin}(n, d/n), \mathrm{Bin}(n-n^{o(1)}, d/n)) \leq 1-(1-d/n)^{n^{o(1)}} = n^{-1+o(1)}$, we can couple these steps exactly with high probability. Applying Bayes' rule, conditioned on $\sigma^\star_{v_i}$, the probability that its neighbor has the same label is $\frac{a}{a+(q-1)b} = \frac{1+(q-1)\lambda}{q}$. Thus, the broadcast process on the tree is chosen so that the labels can be coupled perfectly as well.
\end{proof}

\begin{lemma}\label{lemma: local align}
    There exists $\delta(d) > 0$ such that 
    \[ \Prob{\tau_w = \sigma^\star_w \in \mathcal{C}_\mathrm{good}} > \delta(d). \]
\end{lemma}
\begin{proof}
Let $a, \delta > 0$ be as guaranteed by Lemma~\ref{lemma: branching process} invoked with parameters $n_i$, $d_i$, and $\lambda$. We claim that conditioned on $N_L > 8as_i^L$ and $M_L < as_i^L$ and $\sigma_w^\star \in \mathcal{C}_\mathrm{good}$, we have $\tau_w = \sigma_w^*$ with high probability. Suppose first that $\sigma_w^\star = \sigma_{u_*}^\star \in \mathcal{C}_\mathrm{good}$. With high probability, the number of neighbors of $u_*$ in the same community is at least $\frac{1}{2}\lambda|N(u_*) \cap V_i|$. This guarantees at least $4a\lambda s_i^L|N(u_*)|$ many matches among $N(u_*) \cap V_i$ and $S_i(w)$. 

On the other hand, for any $u_*$ such that $\sigma_w^\star \neq \sigma_{u_*}^\star$, the probability that we have at least $2a\lambda s_i^L|N(u_*)|$ many matches among $N(U_*) \cap V_i$ and $S_i(w)$ is at most the probability that there are at least $2\lambda|N(u_*) \cap V_i|$ neighbors fo $u_*$ with the same community and one of the clusters of $S_i(w)$ being assigned label $\sigma_{u_*}^\star$. This probability is at most 
\[ \exp(-\lambda|N(u_*) \cap V_i|) \cdot \frac{|S_i(w)|}{q} = o\big(\frac{1}{q}\big) \] 
by our choices that $|N(u_*) \cap V_i| \geq \log\log n$ and $L$ growing sufficiently slowly. Union bounding over the remaining representatives gives the desired conclusion.
\end{proof}

The second step shows that the predictions based on path counts tend to align well with the local estimates when they are correct. This is due to the fact that the expected weighted path counts are exactly the local estimators, and the results of Section~\ref{sec: WPS} give tight concentration around this quantity. 
\begin{lemma}\label{lemma: counts align}
    There exists $C$ such that if $\frac{d}{(\log d)^2}\lambda^{1/\chi} > C$, then 
    \[ \Prob{\widehat \sigma_w = \tau_w \middle\vert \sigma^\star_w \in \mathcal{C}_\mathrm{good}} = 1-o(1). \]
\end{lemma}
\begin{proof}
    The first step is to translate the quantities $N_{u,w}^{(i)}$ counting all non-backtracking walks into the quantities $S_{u,w}$ counting only self-avoiding walks. This follows by applying (\ref{eq: complicated paths}) of Proposition~\ref{prop: path count moments}. Notice that among all the $Z_i(v, u_*)$ for $1 \leq i \leq M$, $u_* \in U_*$, and $v \in S_i(w)$, there are a total of at most $q\log n$ terms $N_{u,v}^{(i)}$ with high probability. We need to count how many of these terms are between pairs of vertices $u,v$ in the same community. By Lemma~\ref{lemma: init}, there are at most $\log^2 n$ total such pairs. By (\ref{eq: complicated paths}) and a union bound, each of these $N_{u,v}^{(i)}$ can be replaced with $S_{u,v}^{(i)}$ at a cost of $o(\frac{qs_i^k}{n_i})$ with probability $\log^2 n\left(\frac{q}{n_i}\right)^{1/3} = o(1)$. Union bounding over the remaining pairs, for which the two vertices are in different communities, each $N_{u,v}^{(i)}$ can be replaced by $S_{u,v}^{(i)}$ with probability $q\log^2 n \cdot q^{-2/3}n_i^{-1/3} = o(1)$.  Thus, it suffices to consider the redefined random variables 
    \[ Z_i(w, u_*) \vcentcolon= \sum_{u \in N(u_*)\cap V_i}\sum_{v \in S_i(w)} S_{u,v}^{(i)}, \]
    to which we can apply the moment bounds of Proposition~\ref{prop: path count moments}.

    We control both events for which $\widehat\sigma(w)$ differs from $\tau(w)$. Either $Z_i(w, u_*) - \frac{s_i^k}{n_i}Z_i^\mathrm{loc}(w, u_*)$ is large for some $i$ and $\sigma^\star_w = \sigma^\star_{u_*}$, or $Z_i(w, u) - \frac{s_i^k}{n_i}Z_i^\mathrm{loc}(w, u)$ is large for some $u \neq u_*$ and every $i$. Note that by (\ref{eq: first moment}), $\E{Z_i(w, u_*)\vert \sigma^\star_U} = \frac{s_i^k}{n_i}Z_i^\mathrm{loc}(w, u_*)$, so it suffices to show that $Z_i$ is concentrated around its expectation. This is done via a second moment computation:
    \begin{align*}
        \E{Z_i(w, u_*)^2} &= \E{\left(\sum_{u \in N(u_*)\cap V_i} \sum_{v \in S_i(w)} S_{u,v}^{(i)} \right)^2} = \sum_{u,u'} \sum_{v, v'} \E{S_{u,v}^{(i)}S_{u',v'}^{(i)}}
    \end{align*}
    where both $u$ and $u'$ are in $N(u_*) \cap V_i$ and $v$ and $v'$ are in $S_i(W)$. Thus, we can express
    \[ \Var{Z_i(w, u_*)\middle\vert \sigma^\star_w, \sigma^\star_{u_*}} = \sum_{u,u'}\sum_{v,v'} \E{S_{u,v}^{(i)}S_{u',v'}^{(i)}\middle\vert \sigma^\star_w, \sigma^\star_{u_*}} - \E{S_{u,v}^{(i)}\middle\vert \sigma^\star_w, \sigma^\star_{u_*}}\E{S_{u',v'}^{(i)}\middle\vert \sigma^\star_w, \sigma^\star_{u_*}}. \]

    We partition the pairs in the double summation into groups depending on whether or not $u=u'$, $v=v'$, $\sigma_u^\star = \sigma_v^\star$, and $\sigma_{u'}^\star = \sigma_{v'}^\star$ and control the contribution within each group. 

    \begin{itemize}
        \item $u=u', v=v', \sigma_u^\star = \sigma_v^\star$: By Proposition \ref{prop: path count moments}\eqref{eq: second moment}, contributes $O\big(\frac{q^2s_i^{2k}}{n_i^2}\big)$ to the variance, and there are at most $|N(u_*) \cap V_i||S_i(w)|$ such terms. When $\sigma_w^\star \neq \sigma_{u_*}^\star$, there are at most $O\big(\frac{1}{q}|N(u_*)\cap V_i||S_i(w)|\big)$ such terms in expectation.  
        \item $u=u', v=v', \sigma_u^\star \neq \sigma_v^\star$: By Proposition \ref{prop: path count moments}\eqref{eq: second moment}, contributes $\frac{d_i^k}{n_i}$ to the variance, and there are at most $|N(u_*) \cap V_i||S_i(w)|$ such terms. This will contribute the main term when $\sigma_w^\star \neq \sigma_{u_*}^\star$. 
        \item $u=u', v\neq v', \sigma_u^\star = \sigma_v^\star = \sigma_{v'}^\star$: By Proposition \ref{prop: path count moments}\eqref{eq: second moment} and Cauchy--Schwarz, contributes $O\big(\frac{q^2s_i^{2k}}{n_i^2}\big)$ to the variance, and there are at most $|N(u_*) \cap V_i||S_i(w)|^2$ such terms. When $\sigma_w^\star \neq \sigma_{u_*}^\star$, there are at most $O\big(\frac{1}{q}|N(u_*)\cap V_i||S_i(w)|^2\big)$ such terms in expectation.  
        \item $u=u', v\neq v', \sigma_u^\star = \sigma_v^\star \neq \sigma_{v'}^\star$: By Proposition \ref{prop: path count moments}\eqref{eq: second moment} and Cauchy--Schwarz, contributes $O\big(\frac{d_i^{k/2}qs_i^k}{n_i^{3/2}}\big)$ to the variance, and there are at most $|N(u_*) \cap V_i||S_i(w)|^2$ such terms. When $\sigma_w^\star \neq \sigma_{u_*}^\star$, there are at most $O\big(\frac{1}{q}|N(u_*)\cap V_i||S_i(w)|^2\big)$ such terms in expectation. 
        \item $u=u', v\neq v', \sigma_u^\star \neq \sigma_v^\star, \sigma_{v'}^\star$: By Proposition \ref{prop: path count moments}\eqref{eq: second moment} and Cauchy--Schwarz, contributes $\frac{d_i^k}{n_i}$ to the variance, and there are at most $|N(u_*) \cap V_i||S_i(w)|^2$ such terms.
        \item $u\neq u', v\neq v', \sigma_u^\star = \sigma_{v}^\star, \sigma_{u'}^\star = \sigma_{v'}^\star$: By Proposition \ref{prop: path count moments}\eqref{eq: second moment} and \eqref{eq: cross moment}, contributes $O\big(\frac{q^3s_i^{2k}}{n_i^3}\big)$ to the variance, and there are at most $|N(u_*) \cap V_i|^2|S_i(w)|^2$ such terms. When $\sigma_w^\star \neq \sigma_{u_*}^\star$, there are at most $O\big(\frac{1}{q}|N(u_*)\cap V_i|^2|S_i(w)|^2\big)$ such terms in expectation. 
        \item $u\neq u', v\neq v', \sigma_u^\star = \sigma_{v}^\star, \sigma_{u'}^\star \neq \sigma_{v'}^\star$: By Proposition \ref{prop: path count moments}\eqref{eq: second moment} and \eqref{eq: cross moment}, contributes $O\big(\frac{q^2s_i^{2k}}{n_i^3}\big)$ to the variance, and there are at most $|N(u_*) \cap V_i|^2|S_i(w)|^2$ such terms. 
        \item $u\neq u', v\neq v', \sigma_u^\star \neq \sigma_{v}^\star, \sigma_{u'}^\star \neq \sigma_{v'}^\star$: By Proposition \ref{prop: path count moments}\eqref{eq: second moment} and \eqref{eq: cross moment}, contributes $O\big(\frac{qs_i^{2k}}{n_i^3}\big)$ to the variance, and there are at most $|N(u_*) \cap V_i|^2|S_i(w)|^2$ such terms.
    \end{itemize}
    The cases where $u \neq u', v = v'$ are handled analogously to $u = u', v\neq v'$ with the count $|N(u_*) \cap V_i||S_i(w)|^2$ replaced by $|N(u_*) \cap V_i|^2|S_i(w)|\lambda^L$ in expectation. 

    By the above casework and Chebyshev's inequality we have 
    \begin{multline*}
        \Prob{\abs{Z_i(w, u_*) - \E{Z_i(w, u_*)}} > t\cdot \frac{qs_i^{k}}{n_i}|N(u_*) \cap V_i| \middle\vert \sigma^\star_w = \sigma^\star_{u_*}, S_i(w)} \leq \\
        O\left(\frac{|S_i(w)|^2}{|N(u_*)\cap V_i|t^2} + \frac{|S_i(w)|\lambda^L}{t^2} + \frac{q|S_i(w)|^2}{nt^2} \right)
    \end{multline*}
    in the case when $\sigma_w^\star = \sigma_{u_*}^\star$. We choose $t = \frac{a\lambda s_i^L}{10}$ and for $L$ growing sufficiently slowly the entire right hand side is $o(1)$.

    In the case when $\sigma_w^\star \neq \sigma_{u_*}^\star$, we apply Chebyshev (averaged over conditioning on the labels $\sigma_{S_i(w)}^\star$ and $\sigma_{N(u_*) \cap V_i}^\star$ to deduce
    \begin{multline*}
        \Prob{\abs{Z_i(w, u_*) - \E{Z_i(w, u_*)}} > t\cdot \frac{qs_i^{k}}{n_i}|N(u_*) \cap V_i| \middle\vert \sigma^\star_w \neq \sigma^\star_{u_*}, S_i(w)} \leq \\
        O\left( \frac{|S_i(w)|^2}{q|N(u_*) \cap V_i|t^2} + \frac{d_i^kn_i|S_i(w)|^2}{q^2s_i^{2k}t^2|N(u_*) \cap V_i|} + \frac{d_i^kn_i|S_i(w)|}{q^2s_i^{2k}t^2} + \frac{|S_i(w)|^2}{nt^2} \right) = \frac{d_i^kn_i}{q^2s_i^{2k}} \cdot n^{o(1)}.
    \end{multline*}
    Our value of $k = \lfloor \beta \log n \rfloor$ from Proposition \ref{prop: path count moments} gives the bound $n^{-c_*}M^k$ where 
    \[ c_* = \frac{\chi \frac{\log(\frac{d}{M}\lambda^{1/\chi})}{\log(\frac{d}{M})}}{1 - \chi + \chi\frac{\log(\frac{d}{M}\lambda^{1/\chi})}{\log(\frac{d}{M})}}. \]
    Thus, we conclude that 
    \[ \Prob{\exists i: \abs{Z_i(w, u_*) - \E{Z_i(w, u_*)}} > \frac{1}{10} \frac{aq\lambda s_i^{k+L}}{n_i}|N(u_*) \cap V_i| \middle\vert \sigma^\star_w = \sigma^\star_{u_*}} = o(1) \]
    and 
    \[ \Prob{\abs{Z_i(w, u_*) - \E{Z_i(w, u_*)}} > \frac{1}{10} \frac{aq\lambda s_i^{k+L}}{n_i} |N(u_*) \cap V_i| ,~\forall i \middle\vert \sigma^\star_w \neq \sigma^\star_{u_*}} = n^{-c_*M}M^{kM} = o(\frac{1}{q}). \]
    For this final inequality we require $c_*M - kM\log M > \chi$, which holds for $M \geq C\log d$ and $d\lambda^{1/\chi} > (C\log d)^2$. 
    Taking a union bound over these bad events yields the desired result. 
\end{proof}

\begin{proof}[Proof of Theorem \ref{thm: main alg}(1)]
    Let $m = \lfloor \log n \rfloor$. Let $v_1, \ldots, v_m$ be a uniformly random sample from $V$. By Lemma~\ref{lemma: coupling} and Lemma~\ref{lemma: local align},
    \[ X^\mathrm{loc} \vcentcolon= \frac{1}{m}\sum_{i=1}^{m} \mathbbm{1}\{\tau_{v_i} = \sigma^\star_{v_i} \in \mathcal{C}_\mathrm{good}\} \]
    is concentrated around some $\delta > 0$ with high probability. By Lemma \ref{lemma: counts align}, we know
    \[ X \vcentcolon= \frac{1}{m}\sum_{i=1}^{m} \mathbbm{1}\{\widehat\sigma_{v_i} = \tau_{v_i} = \sigma^\star_{v_i} \in \mathcal{C}_\mathrm{good}\} \]
    satisfies $\E{X} \geq \left(1 - o(1)\right)\E{X^\mathrm{loc}}$. Since $X^\mathrm{loc}-X \geq 0$, Markov's inequality implies 
    \[ \Prob{X^\mathrm{loc}-X > \frac{1}{2} \E{X^\mathrm{loc}}} = o(1). \]
    The concentration of $X^\mathrm{loc}$ implies that 
    \[ \Prob{X^\mathrm{loc}-X > \frac{1}{2}X^\mathrm{loc}} = o(1) \]
    as well, so that $X \geq \frac{1}{2}X^\mathrm{loc} \geq \frac{\delta}{2}$ with probability at least $1 - o(1)$. Finally, since the $\{v_i\}$ are a uniformly random sample, $X$ is concentrated around 
    \[ \frac{1}{n} \sum_{v \in V} \mathbbm{1}\{\widehat\sigma_v = \tau_v = \sigma^\star_v \in \mathcal{C}_\mathrm{good}\} \]
    with high probability. Since this quantity lower bounds $A(\widehat\sigma, \sigma^\star)$, all together this implies that with probability at least $1 - o(1)$, $A(\widehat\sigma, \sigma^\star) \geq \frac{\delta}{2}$.
\end{proof}

\subsection{Above the KS threshold -- Theorem \ref{thm: main alg}(2)}\label{sec: proof thm KS}

\begin{algorithm}[H]\label{alg: above KS}
    \SetAlgoLined
    Take $U \subseteq V$ to be a random subset of size $O(n^{\frac{\chi+1}{2}})$.  \\
     Draw $q$ vertices $U_* = \{u_1, \ldots, u_q\}$ uniformly at random from all vertices in $U$ having degree at least $\log\log n$ to $V\setminus U$. \\
    \For{$w \in V \setminus U$}{
    \For{$\ell \in [q]$}{
    Compute $Z(w, u_\ell) = \sum_{u \in N(u_\ell)} \sum_{v \in S(w)} N_{u,v}$. \\
    }
    Set $\widehat{\sigma}_w = \ell$ if $\ell$ is the unique label such that $Z(w, u_\ell) > 3aq\lambda \frac{s^{k+L}}{n}|N(u_\ell)| - \frac{s^k}{n}|N(u_\ell)||S(w)|$.
    }
    \For{ $w \in U \setminus U_*$ } {
    Assign $\widehat\sigma(w)$ uniformly at random.
    }
    \For{$\ell \in [q]$}{
    Assign $\widehat\sigma(u_\ell) = \ell$.
    }
    \caption{Above the KS threshold}
\end{algorithm}

In this subsection we analyze Algorithm~\ref{alg: above KS} and prove its correctness. In the regime above the KS threshold, the non-backtracking walk statistics are sufficiently concentrated to ensure small candidate lists (in fact unique candidates) for most vertices without the additional boosting via partitioning the graph. Thus, we can bypass the boosting step and directly compute using the entire graph. In particular, we compute 
\[ Z(w, u_*) = \sum_{u \in N(u_*)}\sum_{v \in S(w)} N_{u,v} \]
and assign $u_*$ as a candidate for $w$ if this value is large. 

As before let $\mathcal{C}_\mathrm{good}$ be the set of good communities that are uniquely represented and do not appear too often in each other's neighborhoods. We once again compare with the local estimator to show that the candidate vertices align well with $\sigma^\star$. Notice that Claim~\ref{claim} (without partitioning) and Lemmas~\ref{lemma: local align} and~\ref{lemma: coupling} remain true, so it suffices to prove an analog of Lemma~\ref{lemma: counts align}.
\begin{lemma}\label{lem: counts align KS}
    If $d\lambda^2 > 1$, then $\Prob{\widehat\sigma_v = \tau_v | \sigma^\star_v \in G} = 1-o(1)$.
\end{lemma}
\begin{proof}
    The proof is identical to the proof of Lemma~\ref{lemma: counts align} excluding the computation of the second moment. We perform the new second moment computation here, and refer the reader to the previous proof for the remainder of the argument. 
    
    We compute the second moment with the same casework, using Proposition \ref{prop: path count moments KS} rather than \ref{prop: path count moments}\eqref{eq: second moment}. In the setting of $\sigma_{u_*}^\star = \sigma_w^\star$, the main contributors are the same (as $\frac{qs^{2k}}{n^2}$ is $o(\frac{q^2s^{2k}}{n^2})$ just like $\frac{d^k}{n}$ was), so the proof proceeds identically and we obtain the same bound of $o(1)$. 

    In the setting of $\sigma_{u_*}^\star \neq \sigma_w^\star$, the main contribution from the term $\frac{d^k}{n}$ is replaced by a contribution of $\frac{qs^{2k}}{n^2}$. Thus, the corresponding bound becomes
    \begin{multline*}
        \Prob{\abs{Z(w, u_*) - \E{Z(w, u_*)}} > t\cdot \frac{qs^{k}}{n}|N(u_*)| \middle\vert \sigma^\star_w \neq \sigma^\star_{u_*}} \leq \\
        O\left( \frac{|S(w)|^2}{q|N(u_*)|t^2} + \frac{|S(w)|^2}{qt^2|N(u_*)|} + \frac{|S_i(w)|}{qt^2} + \frac{|S(w)|^2}{nt^2} \right) = o\big(\frac{1}{q}\big)
    \end{multline*}
    as $\E{|S(w)|} = d^L$ and we chose $t^2 = \frac{a^2\lambda^2}{100}s^{2L} \gg d^L$ when $L$ is (slowly) growing with $n$. Thus, the union bound to conclude the proof follows immediately.
\end{proof}

\subsection{Low-degree implications -- Corollary \ref{cor: sparse transition}(1)}\label{subsec: LDI}
In this subsection we interpret our algorithmic achievability results in the low-degree framework. We prove that the self-avoiding walk statistic is both low degree and has constant correlation with the true indicators. This follows from a direct computation given the second moment estimates in Subsection~\ref{sec: WPS}.
\begin{proof}[Proof of Corollary \ref{cor: sparse transition}(1)]
By Proposition~\ref{prop: path count moments KS} we have the estimates
\begin{align*}
    \Cov{S_{u,v}}{\mathbbm{1}\{\sigma^\star_u = \sigma^\star_v\} - \frac{1}{q}} &= \frac{1}{q}\cdot\left(1-\frac{1}{q}\right)\E{S_{u,v} \vert \sigma^\star_u=\sigma^\star_v} + \frac{q-1}{q}\cdot\left(-\frac{1}{q}\right)\E{S_{u,v} \vert \sigma^\star_u\neq\sigma^\star_v} \\
    &= (1+o(1))\frac{1}{q}\cdot\left(1-\frac{1}{q}\right)\left( \frac{(q-1)s^k}{n} + \frac{s^k}{n} \right) = (1+o(1))\frac{q-1}{q} \cdot \frac{s^k}{n} \\
    \Var{S_{u,v}} &= \E{S_{u,v}^2} \leq \frac{1}{q} \cdot \frac{s^2}{(s-1)^2}\cdot \frac{q^2s^{2k}}{n^2} + \frac{q-1}{q} \cdot \frac{s^2}{s^2-d}\cdot \frac{qs^{2k}}{n^2} \\
    &= \left(\frac{s^2}{(s-1)^2} + \frac{s^2}{s^2-d}\right)\frac{qs^{2k}}{n^2} \\
    \Var{\mathbbm{1}\{\sigma^\star_u=\sigma^\star_v\} - \frac{1}{q}} &= \frac{1}{q}\left(1-\frac{1}{q}\right).
\end{align*}
Notice that $S(u,v)$ is a degree-$\beta\log n$ estimator where $\beta$ is chosen as in Proposition~\ref{prop: path count moments KS}. All together, we have that for $\chi < \frac{1}{2}$ and $d\lambda^{2} > 1$, 
\begin{align*}
    \Corr_{\leq \beta\log n} &\geq (1+o(1))\frac{\frac{s^k}{n}}{\sqrt{\left(\frac{s^2}{(s-1)^2} + \frac{s^2}{s^2-d}\right)\frac{qs^{2k}}{n^2} \cdot \frac{1}{q}}} \geq \frac{1+o(1)}{\sqrt{\frac{s^2}{(s-1)^2} + \frac{s^2}{s^2-d}}} = \delta(d, s).
\end{align*}
\end{proof}
We remark that a similar computation shows that when $\chi > \frac{1}{2}$, $d\lambda^{1/\chi} > C(\log d)^2$, the degree-$\Theta(\log d\log n)$ correlation is also bounded below by a positive constant independent of $n$. This could serve as a partial converse to the low-degree hardness result in this regime.

\section{Low-degree hardness}\label{sec: LDH}

In this section, we prove Theorem~\ref{thm: main hard} and Theorem~\ref{thm:spase:mainhard}. Throughout this section, we consider both the sparse and the dense regimes. To deal with the latter case, we consider the \emph{adjusted} average degree and signal strength:
\[
d_{\circ}:=\frac{a(1-\frac{a}{n})+(q-1)b(1-\frac{b}{n})}{q},\quad\quad \la_{\circ}:=\frac{a-b}{a(1-\frac{a}{n})+(q-1)b(1-\frac{b}{n})}.
\]
Here, we note that the difference between $d_{\circ}$ (resp.\ $\la_{\circ}$) and $d\equiv \frac{a+(q-1)b}{q}$ (resp.\ $\la \equiv \frac{a-b}{a+(q-1)b}$) becomes relevant only when $a\asymp n$ or $b\asymp n$. Define the parameter
\begin{equation}\label{eq:def:xi}
\xi:=\frac{\min\{a(1-\frac{a}{n}),b(1-\frac{b}{n})\}}{\dc}.
\end{equation}
Theorem~\ref{thm: main hard} and Theorem~\ref{thm:spase:mainhard} are consequences of the following result.
\begin{theorem}\label{thm:sbm:low:degree}
Consider the degree-$D$ correlation $\Corr_{\leq D}$ defined in Eq.~\eqref{eq:def:corr}. There exist universal constants $c,C>0$ such that the following holds.
\begin{enumerate}
\item If
$D\leq \left(\frac{\xi n}{Cq^2} \right)^c$, then
\[
\Corr_{\leq D}^2\le \frac{Cq}{n}\sum_{\ell=1}^{D}(\dc \la_{\circ}^2)^{\ell}.
\]
\item If $D\leq \left(\frac{n}{Cq}\right)^c$, then
\[
\Corr_{\leq D}^2\leq \frac{Cq}{n}\sum_{\ell=1}^{D}\left(1+\frac{(q-1)b(1-\frac{b}{n})}{a(1-\frac{a}{n})}\right)^{\ell}\cdot \left(d_{\circ}\la_{\circ}^2\right)^{\ell}.
\]
\end{enumerate}
\end{theorem}
We first show that Theorem~\ref{thm:sbm:low:degree} implies Theorem~\ref{thm: main hard} and Theorem~\ref{thm:spase:mainhard}.
\begin{proof}[Proof of Theorem~\ref{thm: main hard}]
For simplicity, let us denote $a_{\circ}=a(1-a/n)$ and $b_{\circ}=b(1-b/n)$. By the second item of Theorem~\ref{thm:sbm:low:degree}, if $d_{\circ}\la_{\circ}^2 (1+\frac{(q-1)b_{\circ}}{a_{\circ}})\leq 1-c'$ for some $c'=c'(\eta)$, then the desired conclusion holds. Note that if $(q-1)b_{\circ}\leq \eta a_{\circ}$, this condition is clearly satisfied with $c'=\eta^2$ since $\text{SNR}=d_{\circ}\la_{\circ}^2\leq 1-\eta$. As a consequence, we are left to consider the case $(q-1)b_{\circ}\geq \eta a_{\circ}$. In this case, note that
\[
\xi\equiv \frac{q\cdot(a_{\circ}\wedge b_{\circ})}{a_{\circ}+(q-1)b_{\circ}}\geq \frac{q}{(q-1)\eta^{-1}+q-1}\left(\frac{a_\circ}{b_{\circ}}\wedge 1\right).
\]
Also, because we assumed that $a\geq \delta b$ and $a\leq n/2$, we have $a_{\circ}\geq \frac{a}{2}\geq \frac{\delta b_{\circ}}{2}$. Thus, it follows that
\[
\xi \geq \frac{q}{q-1}\cdot \frac{\frac{\delta}{2} \wedge 1}{\eta^{-1}+1}\geq \frac{\frac{\delta}{2} \wedge 1}{\eta^{-1}+1},
\]
so the first item of Theorem~\ref{thm:sbm:low:degree} concludes the proof.
\end{proof}
\begin{proof}[Proof of Theorem~\ref{thm:spase:mainhard}]

Suppose $q=\omega(1), a=\Theta(q)$, and $b=\Theta(1)$. Then, $d_{\circ}=d\left(1+O\left(\frac{q}{n}\right)\right)$ and $\la_{\circ}=\la \left(1+O\left(\frac{q}{n}\right)\right)$, so $d\la^2\leq 1$ implies that $d_{\circ}\la_{\circ}^2\leq 1+O\left(\frac{q}{n}\right)$. Moreover, since $\min\{a(1-\frac{a}{n}),b(1-\frac{b}{n})\}=\Theta(1)$ and $d_{\circ}=\Theta(1)$, we necessarily have $\xi=\Omega(1)$. Thus, the desired claim follows from the first item of Theorem~\ref{thm:sbm:low:degree}.
\end{proof}
For the rest of this section, we prove Theorem~\ref{thm:sbm:low:degree}. Following the framework of~\cite{SW:24}, our proof will take the following steps. We will choose a basis $\{\phi_\alpha\}_{\alpha \in \cI}$ for $\R[Y]_{\le D}$, so that an arbitrary degree-$D$ polynomial (in $Y$) can be expanded as
\[ f(Y) = \sum_\alpha \hat{f}_\alpha \phi_\alpha(Y). \]
We will choose a collection of random variables $\{\psi_{\beta\gamma}\}_{(\beta,\gamma) \in \cJ}$ that are functions of $(Y,\sigma^\star)$, and are furthermore orthonormal in the sense $\EE[\psi_{\beta\gamma} \cdot \psi_{\beta'\gamma'}] = \one_{\beta=\beta'} \cdot \one_{\gamma=\gamma'}$. This allows us to bound $\EE[f(Y)^2] \ge \sum_{\beta\gamma} \EE[f(Y) \cdot \psi_{\beta\gamma}]^2 = \|M \hat f\|^2$ where
\[ M_{\beta\gamma,\alpha} := \EE[\phi_\alpha(Y) \cdot \psi_{\beta\gamma}]. \]
Also define
\[ c_\alpha := \EE[\phi_\alpha(Y) \cdot x]~~~\textnormal{where}~~~ x \equiv \one\{\sigma^\star_1=\sigma^\star_2\}-\frac{1}{q}.
\]
We will construct a vector $u = (u_{\beta\gamma})_{(\beta,\gamma) \in \cJ}$ such that $u^\top M = c^\top$, which yields the bound
\begin{equation}\label{eq:corr-bound}
\sqrt{\EE[x^2]} \cdot \Corr_{\le D} = \sup_f \frac{\EE[f(Y) \cdot x]}{\sqrt{\EE[f(Y)^2]}} \leq \sup_{\hat f} \frac{c^\top \hat f}{\|M \hat f\|} = \sup_{\hat f} \frac{u^\top M \hat f}{\|M \hat f\|} \leq \sup_{\hat f} \frac{\|u\| \cdot \|M \hat f\|}{\|M \hat f\|} = \|u\|.
\end{equation}
The constraints $u^\top M = c^\top$ can be equivalently expressed as
\begin{equation}\label{eq:u-cond}
\sum_\beta u_\beta M_{\beta\gamma,\alpha} = c_\alpha.
\end{equation}

\subsection{Choosing the basis}
For $\alpha = (\alpha_{ij})_{1 \le i < j \le n} \in \{0,1\}^{\binom{n}{2}}$, we consider the basis $\{\phi_{\al}:\al\in \{0,1\}^{\binom{n}{2}}\}$ defined as
\[
\phi_{\al}(Y):=\prod_{1 \le i < j \le n}\left(Y_{ij}-\frac{b}{n}\right)^{\al_{ij}}.
\]
We remark that this choice differs from~\cite{SW:24}, where they instead considered $(Y-\dc/n)^{\al}$. As shown in Eq.~\eqref{eq:compute:M} below, this new basis facilitates an explicit calculation of $M_{\be\ga,\al}$, thereby yielding a sharper upper bound on $|M_{\be\ga,\al}|$ (see Lemma~\ref{lem:bound:M}). We note that~\cite{SW:24} focused on the case of constant $q$ but treated more general block models with communities of varying sizes and varying connection probabilities (and in that setting, there is not a clear analogue of our new basis, since the across-community edges do not all have the same probability $b/n$). Our focus here is instead on the symmetric block model with growing $q$. 

The choice for the orthonormal random variables will be the same as~\cite{SW:24}. Note that $\be\in \{0,1\}^{\binom{n}{2}}$ can be viewed as a graph on vertex set $[n]$, namely the graph containing each edge $(i,j)$ for which $i<j$ and $\beta_{ij} = 1$. Let $E(\beta)$ denote the set of these edges $(i,j)$, and let $V(\be) \subseteq [n]$ denote the set of non-isolated vertices of this graph. For $\be \in \{0,1\}^{\binom{n}{2}}$ and $\ga\in [q]^{V(\be)}$, let
\[
\psi_{\be\ga}(Y,\sigma^\star):=\one\big\{\sigma^\star_v=\gamma_v,~\forall v\in V(\be)\big\} \, q^{|V(\beta)|/2} \prod_{\substack{(i,j)\in E(\be)\\\gamma_i= \gamma_j}}\bigg\{\frac{Y_{ij}-\frac{a}{n}}{\sqrt{\frac{a}{n}(1-\frac{a}{n})}}\bigg\} \prod_{\substack{(i,j)\in E(\be)\\\gamma_i\neq  \gamma_j}}\bigg\{\frac{Y_{ij}-\frac{b}{n}}{\sqrt{\frac{b}{n}(1-\frac{b}{n})}}\bigg\}.
\]
If $\beta=\emptyset$, we take the convention $\ga=\emptyset$. The following lemma follows from \cite{SW:24}.

\begin{lemma}[Lemma 6.2 in \cite{SW:24}]\label{lem:psi:ortho}
    The set $\{\psi_{\be\ga}: \be \in \{0,1\}^{\binom{n}{2}}, \ga\in [q]^{V(\be)}\}$ is orthonormal.
\end{lemma}

An important consequence of our choice of $\{\phi_{\al}\}$ and $\{\psi_{\be\ga}\}$ is the following. Recall
\[
M_{\be\ga,\al} \equiv \EE \big[\phi_{\al}(Y)\psi_{\be\ga}(Y,\sigma^\star)\big],
\]
and for $1 \le i < j\le n$, define
  \begin{equation}\label{eq:def:Z}
    Z_{ij}:=
    \begin{cases}
        \frac{Y_{ij}-\frac{a}{n}}{\sqrt{\frac{a}{n}(1-\frac{a}{n})}} &~~~\textnormal{if}~~\sigma^\star_i=\sigma^\star_j,\\
         \frac{Y_{ij}-\frac{b}{n}}{\sqrt{\frac{b}{n}(1-\frac{b}{n})}}&~~~\textnormal{otherwise.}
    \end{cases}
    \end{equation}
\begin{lemma}\label{lem:M}
$M_{\be\ga,\al}=0$ unless $\be \leq \al$.
\end{lemma}
Here $\beta \le \alpha$ is defined entrywise, i.e., it means $\beta_{ij} \le \alpha_{ij}$ for all $i < j$.

\begin{proof}
Suppose $(u,v) \in E(\beta)\setminus E(\al)$. Recalling $Z_{ij}$ from \eqref{eq:def:Z}, let $\psi_{\be\ga}^{-}:=\psi_{\be\ga}/Z_{uv}$ and write $Y^{-}:=(Y_{ij})_{\{i,j\}\neq \{u,v\}}$. Then, note that $\phi_{\al}$ and $\psi^{-}_{\be\ga}$ are $(Y^{-},\sigma^\star)$-measurable. Since $\EE[Z_{ij}\given \sigma^\star_u,\sigma^\star_v]=0$, 
\[
\EE[\phi_{\al}\psi_{\beta\gamma}\given Y^{-},\sigma^\star]=\phi_{\al}\psi^{-}_{\be\ga}\cdot \EE\big[Z_{uv} \bgiven \sigma^\star_u,\sigma^\star_v\big]=0,
\]
so taking expectations on both sides shows $M_{\be\ga,\al}=0$.
\end{proof}

\subsection{Removing the uninformative terms}

\begin{definition}\label{def:informative}
         Define the set $\hat\cI\subseteq \{0,1\}^{\binom{n}{2}}$ of ``informative'' $\alpha$ as the ones satisfying the following conditions, with $\emptyset \in \hat\cI$ by convention:
    \begin{itemize}
        \item $1,2\in V(\al)$,
        \item for all $v \in V(\al) \setminus \{1,2\}$, $\deg_{\al}(v)\geq 2$,
        \item $\al$ is connected.
    \end{itemize}
    Also define $\hat\cJ:=\{(\be,\ga) \,:\, \be\in \hat\cI,\, \ga\in [q]^{V(\be)}\}$ with the convention that if $\beta=\emptyset$ then $\gamma=\emptyset$.
\end{definition}

\begin{lemma}\label{lem:colinearity}
    For each $\al \notin \hat\cI$ there exists $\hat\al\in \hat\cI$ and $\mu \in \R$ such that $c_{\al}=\mu c_{\hat\al}$ and $M_{\be\ga,\al}=\mu M_{\be\ga,\hat\al}$ for all $(\be,\ga) \in \hat\cJ$.
\end{lemma}
\noindent The purpose of this lemma is as follows. We will set $u_{\beta\gamma} = 0$ for all $(\beta,\gamma) \notin \hat\cJ$. Then Lemma~1.4 of~\cite{SW:24} implies that we only need to verify the constraints~\eqref{eq:u-cond} for $\alpha \in \hat\cI$, as the constraints for $\alpha \notin \hat\cI$ then follow automatically. As a result, it suffices to verify~\eqref{eq:new-u-cond} below.

\begin{proof}
Consider $\al\notin \hat\cI$. We divide into two cases. First, suppose that $\{1,2\}\nsubseteq V(\al)$, and w.l.o.g.\ assume $2\notin V(\al)$. By setting $\hat \al=\emptyset$, Lemma~\ref{lem:M} shows that $M_{\be\ga,\al}=0=M_{\be\ga, \hat\al}$ for all $\be \neq \emptyset$. Moreover, since $2\notin V(\al)$, 
\[
c_{\al}=\EE \Big[\EE\big[\phi_{\al}(Y)\bgiven \sigma^\star_1\big]\cdot \EE\Big[\one\{\sigma^\star_1=\sigma^\star_2\}-\frac{1}{q}\Bgiven \sigma^\star_1\Big]\Big]=0,
\]
where the last equality holds because $\P(\sigma^\star_1=\sigma^\star_2\given \sigma^\star_1)=q^{-1}$ a.s.\ by symmetry.  Thus, $c_{\al}=\mu c_{\emptyset}$ and $M_{\be\ga,\al}=\mu M_{\be\ga,\emptyset}$ for $\be\in \hat\cI, \ga\in [q]^{V(\be)}$ are satisfied by setting $\mu=M_{\emptyset,\al}/M_{\emptyset,\emptyset}$.

Next, suppose that $1,2\in V(\al)$ but either $\al$ is not connected or there exists $v\in V(\al)\setminus\{1,2\}$ such that $\deg_{\al}(v)=1$. We then consider $\hat\al$ produced by the following steps.
\begin{itemize}
    \item At time $t=0$, let $\al_0$ be the connected component of $\al$ that contains both of the vertices $1,2$. If none, we set $\hat\al=\emptyset$ and terminate.
    \item At time $t\geq 1$, if there exists a vertex $v_{t-1}\in V(\al_{t-1})\setminus\{1,2\}$ such that $\deg_{\al_{t-1}}(v_{t-1})=1$, then set $\al_{t}=\al_{t-1}-\bone_{(u_{t-1}v_{t-1})}$ where $(u_{t-1},v_{t-1})\in E(\al_{t-1})$ is the edge adjacent to $v_{t-1}$. That is, at time $t$, trim a leaf of $\al_{t-1}$, if any.
    \item Iterate the previous step until done.
\end{itemize}
Let $\hat \al$ be the resulting graph which must be informative (Definition~\ref{def:informative}). It suffices to show that there exists $\mu_t\in \R, t\geq -1$ such that $c_{\al_{t}}=\mu_t c_{\al_{t+1}}$ and $M_{\be\ga, \al_t}=\mu_t M_{\be\ga, \al_{t+1}}$ holds for any $\be\in \hat\cI$ and $\ga\in [q]^{V(\be)}$, where $\al_{-1}\equiv \al$ by convention. Here, note that $\al_{t+1}\leq \al_{t}$ holds by construction and so by Lemma~\ref{lem:M}, if $\be\nleq \al_{t}$, then $M_{\be\ga,\al_{t}}=M_{\be\ga,\al_{t+1}}=0$ holds. To this end, we restrict our attention to $\be\leq \al_{t}$.

At time $t=0$, consider $\be\in \hat\cI$ such that $\be\le \al$. Since a non-empty $\be\in \hat\cI$ is connected with $1,2\in V(\be)$, we always have $\be\leq \al_0$. Moreoever, since $V(\al_0)\cap V(\al-\al_0)=\emptyset$, we have
\[
M_{\be\ga, \al}=\EE[\phi_{\al-\al_0}\cdot \psi_{\al_0}\psi_{\be\ga}]=\EE[\phi_{\al-\al_0}]\cdot \EE[\psi_{\al_0}\psi_{\be\ga}]\equiv \mu_0\cdot M_{\be\ga,\al_0}.
\]
Similarly, in the case where $1,2\in V(\al_0)$, $c_{\al}=\EE[\phi_{\al-\al_0}]\cdot \EE[\phi_{\al_0} x]= \mu_0\cdot c_{\al_0}$ holds by the independence of components. In the other case, the vertices $1$ and $2$ are contained in different connected components of $\al$. Note that by symmetry of the colors, revealing spin labels of one vertex in each connected component of $\al$ does not affect the expected value of $\phi_{\al}$. In particular, we have $\EE[\phi_{\al}\given \sigma^\star_1=\sigma^\star_2]= \EE[\phi_{\al}]$ in this case, so $c_{\al}= \EE[\phi_{\al}]\cdot \EE[x]= \mu_0\cdot c_{\emptyset}$. Hence, in either case, $c_{\al}=\mu_0\cdot c_{\al_0}$ holds.

At time $t\geq 1$, let $\be\in \hat\cI$ such that $\be \leq \al_t$. Such $\be$ cannot contain the degree-$1$ vertex $v_t\in V(\al_t)\setminus \{1,2\}$, so $\be \leq \al_{t-1}$. Note that revealing the spin label of the vertex $u_t$ adjacent to $v_t$ makes $\phi_{\al_{t+1}}\psi_{\be\ga}$ and $Y_{u_t v_t}$ conditionally independent. Thus,
\[
M_{\be\ga,\al_t}=\EE\Big[\EE\big[\phi_{\al_{t+1}}\psi_{\be\ga}\bgiven \sigma^\star_{u_t}\big]\cdot \EE\big[Y_{u_t v_t}-b/n\bgiven \sigma^\star_{u_t}\big]\Big].
\]
An important observation is that by symmetry of the colors, $\phi_{\al_{t+1}}\psi_{\be\ga}$ is independent from $\sigma^\star_{u_t}$, thus $\EE\big[\phi_{\al_{t+1}}\psi_{\be\ga}\bgiven \sigma^\star_{u_t}\big]=M_{\be\ga,\al_{t+1}}$, a.s. Hence,
\[
M_{\be\ga,\al_t}= M_{\be\ga,\al_{t+1}}\cdot \EE\big[Y_{u_t v_t}-b/n\big]=M_{\be\ga,\al_{t+1}}\cdot \frac{a-b}{nq}.
\]
By the same argument,
\[
c_{\al_t}= \EE\Big[\EE\big[\phi_{\al_{t+1}}x\bgiven \sigma^\star_{u_t}\big]\cdot \EE\big[Y_{u_t v_t}-b/n\bgiven \sigma^\star_{u_t}\big]\Big]=c_{\al_{t+1}}\cdot \frac{a-b}{nq},
\]
which concludes the proof.
\end{proof}

\subsection{Construction of \texorpdfstring{$u$}{u}}
By Lemma~\ref{lem:colinearity} (see also Lemma~1.4 of~\cite{SW:24}), our goal is now to construct $u=(u_{\be\ga})_{\be\ga\in \hat\cJ}$ such that
\begin{equation}\label{eq:new-u-cond}
    \sum_{\be\ga \in \hat\cJ} u_{\be\ga}M_{\be\ga,\al}=c_{\al}\qquad\forall\al\in \hat\cI,
\end{equation}
with $\|u\|$ as small as possible to obtain a good bound on $\Corr_{\le D}$. We will write 
\[
|\al| := \sum_{i<j} \alpha_{ij} = |E(\alpha)|.
\]
We construct $u$ in the same manner as in \cite[Section 6.3]{SW:24}. Namely, it is constructed implicitly by the following recursion:
\begin{enumerate}
    \item Set $u_{\emptyset}\equiv c_{\emptyset}=0$. 
    \item Having  determined $u_{\be\ga}$ for all $\be\in \hat\cI$ such that $|\be|<|\al|$, set
    \begin{equation}\label{eq:def:d}
   d_{\al}:=c_{\al}-\sum_{\be\ga\in \hat\cJ \,:\, \be\lneq\al}u_{\be\ga}M_{\be\ga,\al}
    \end{equation}
    and
\begin{equation}\label{eq:sbm:u}
    u_{\al\ga}=\frac{M_{\al\ga,\al}}{\sum_{\ga'\in [q]^{V(\al)}}M_{\al \ga', \al}^2}d_{\al}, \qquad \ga\in [q]^{V(\al)}.
    \end{equation}
\end{enumerate}
The intuition behind this choice of $u_{\al\ga}$ is that given $\{u_{\be\ga} \,:\, \be\ga\in \hat\cJ, \, |\be|<|\al|\}$, the choice~\eqref{eq:sbm:u} minimizes $\sum_{\ga\in [q]^{V(\al)}}u_{\al\ga}^2$ by Cauchy--Schwarz. Recalling~\eqref{eq:corr-bound}, we conclude
\begin{equation}\label{eq:sbm:corr:bound}
\Corr^2_{\leq D}\leq \frac{\|u\|^2}{\EE[x^2]}=\frac{q}{1-1/q} \sum_{\al\in \hat\cI \, : \, 1\leq |\al|\leq D} \frac{d_\al^2}{\sum_{\ga\in [q]^{V(\al)}}M_{\al\ga,\al}^2}.
\end{equation}
To complete the proof, we next upper bound $d_{\al}^2$ and lower bound $\sum_{\ga\in [q]^{V(\al)}}M_{\al\ga,\al}^2$.

We first compute $c_{\al}\equiv \EE[\phi_{\al}(Y)\cdot x]$ for $\al\in \hat\cI$. Note that $\EE[Y_{ij}-\frac{b}{n}\given \sigma^\star]=\frac{a-b}{n}\cdot \one\{\sigma^\star_i=\sigma^\star_j\}$, so
\[
\EE\big[\phi_{\al}(Y)\cdot x \bgiven \sigma^\star\big]=\left(\frac{a-b}{n}\right)^{|\al|}\cdot \left(\one\{\sigma^\star_1=\sigma^\star_2\}-\frac{1}{q}\right)\cdot \one\{\sigma^\star_i=\sigma^\star_j,~~\forall i,j\in V(\al)\}.
\]
Since $\al\in \hat\cI$ is connected, it follows that
\begin{equation}\label{eq:compute:c}
    c_{\al}=(q-1)\cdot q^{-|V(\al)|}\cdot \left(\frac{a-b}{n}\right)^{|\al|}.
\end{equation}
Next, we compute $M_{\be\ga,\al}$ for $\be \leq \al$. Recalling the notation $Z_{ij}$ from Eq.~\eqref{eq:def:Z}, note that
\[
\EE\big[Z_{ij}\cdot (Y_{ij}-b/n)\bgiven \sigma^\star\big]=
\begin{cases}
   \sqrt{\frac{a}{n}\left(1-\frac{a}{n}\right)}&\quad\text{if}\quad \sigma^\star_i=\sigma^\star_j,\\
   \sqrt{\frac{b}{n}\left(1-\frac{b}{n}\right)}&\quad\text{if}\quad \sigma^\star_i\neq \sigma^\star_j.
\end{cases}
\]
Thus, for $\be\le \al$, a direct calculation shows
\begin{equation*}
\begin{split}
\EE\big[\phi_{\al}\psi_{\be\ga}\bgiven \sigma^\star\big]
&=q^{\frac{|V(\be)|}{2}}\left(\frac{a}{n}\left(1-\frac{a}{n}\right)\right)^{\frac{\ell(\be,\ga)}{2}}\left(\frac{b}{n}\left(1-\frac{b}{n}\right)\right)^{\frac{|\be|-\ell(\be,\ga)}{2}} \left(\frac{a-b}{n}\right)^{|\al|-|\be|}\\
&\qquad\qquad\cdot \one\Big\{\sigma^\star_v=\gamma_v,~\forall v\in V(\be)~~~\textnormal{and}~~~~\sigma^\star_u=\sigma^\star_v,~\forall(u,v)\in E(\al)\setminus E(\be)\Big\},
\end{split}
\end{equation*}
where $\ell(\be,\ga)$ denotes the number of monochromatic edges in $\be$ (i.e., edges in $\be$ belonging to the same community) under the community labeling $\ga\in [q]^{V(\be)}$:
\begin{equation*}
    \ell(\be,\ga)=\big|\big\{(u,v)\in E(\beta) : \ga_u=\ga_v\big\}\big|.
\end{equation*}
We can thus express $M_{\be\ga,\al}$ as
\begin{equation}\label{eq:compute:M}
\begin{split}
    M_{\be\ga,\al}
    &=q^{-\frac{|V(\be)|}{2}}\left(\frac{a}{n}\left(1-\frac{a}{n}\right)\right)^{\frac{\ell(\be,\ga)}{2}}\left(\frac{b}{n}\left(1-\frac{b}{n}\right)\right)^{\frac{|\be|-\ell(\be,\ga)}{2}} \left(\frac{a-b}{n}\right)^{|\al|-|\be|}\\
    &\qquad\qquad\quad \cdot \P\Big(\sigma^\star_u=\sigma^\star_v,~\forall(u,v)\in E(\al)\setminus E(\be)\Bgiven \sigma^\star_v=\gamma_v,~\forall v\in V(\be)\Big).
\end{split}
\end{equation}

\noindent With the expression~\eqref{eq:compute:M} in hand, we prove the following estimate, which plays a crucial role in recursively upper bounding $|d_{\al}|$.
\begin{lemma}\label{lem:bound:M}
    For any $\al,\be\in \hat\cI$ such that $\be\le \al$ and $\ga\in [q]^{V(\be)}$, we have
    \[
    \left|M_{\be\ga,\al}\right|\leq M_{\be\ga,\be}\cdot\left(\frac{|a-b|}{n}\right)^{|\al|-|\be|} q^{-|V(\al)|+|V(\be)|}.
    \]
\end{lemma}
\begin{proof}
Immediately from the expression~\eqref{eq:compute:M}, we have $M_{\be\ga,\be}\geq 0$. Moreover, 
   \[
   \left|M_{\be\ga,\al}\right|=M_{\be\ga,\be}\cdot \left(\frac{|a-b|}{n}\right)^{|\al|-|\be|}\P\Big(\sigma^\star_u=\sigma^\star_v,~\forall(u,v)\in E(\al)\setminus E(\be)\Bgiven \sigma^\star_v=\gamma_v,~\forall v\in V(\be)\Big).
   \]
   Observe that the conditional probability above is either $q^{-(|V(\al)|-|V(\be)|)}$ or $0$. Indeed, let $\al_1,\ldots, \al_K$ be the connected components of $\al-\be$. Then, the conditional probability is $q^{-(|V(\al)|-|V(\be)|)}$ if for all $i\in [K]$, $\gamma_u=\gamma_v$ holds for any $u,v \in V(\al_i)$, and $0$ otherwise. This concludes the proof.
\end{proof}

The following recursion $f(\cdot)$ plays a crucial role in upper bounding $|d_{\al}|$: observe that the informative $\al\in \{0,1\}^{\binom{\N}{2}}$ can be defined analogously by considering $n=\infty$. Let $\hat\cI_{\infty}$ denote the corresponding set of informative $\al$'s, and let $f:\hat\cI_{\infty}\to \N$ be defined by the recursion
\begin{equation}\label{eq:def:recursion}
f(\alpha)=\sum_{\substack{0 \le \beta\lneq\alpha \\ \beta\in \hat\cI_{\infty}}} f(\beta) \qquad \text{for } \alpha \ne \emptyset,
\end{equation}
with the base case $f(\emptyset)=1$. Intuitively, $f(\al)$ measures the ``complexity'' of $\al$, and the following lemma was proved in~\cite{SW:24} that $f(\al)$ grows at most exponentially in the number of tree excess edges.
\begin{lemma} [Lemma 5.4 in \cite{SW:24}] \label{lem:f:bound}
We have $f(\alpha)\leq (2 |\alpha|)^{|\alpha|-|V(\alpha)|+1}$ for any non-empty $\al\in \hat\cI_{\infty}$.
\end{lemma}
\begin{remark}\label{rmk:comparison:sw}
    Note that $\hat\cI$, the set of ``informative'' $\al$, is slightly different from the set of ``good'' $\al$ in~\cite{SW:24} (see Definition~6.3 therein). Namely, our ``informative'' $\al$ requires a slightly stronger condition that $\al$ is connected compared to $\al+\mathbf{1}_{(1,2)}$ is connected, where $\mathbf{1}_{(1,2)}$ is the graph consisting of a single edge $(1,2)$. Since the recursion~\eqref{eq:def:recursion} leads to a larger value of $f(\cdot)$ when $\hat\cI_{\infty}$ is enlarged, Lemma~\ref{lem:f:bound} is indeed a consequence of Lemma~5.4 in~\cite{SW:24}.
\end{remark}
Next, we establish an upper bound on $|d_{\al}|$ depending on $f(\al)$.

\begin{lemma}\label{lem:d:alpha}
    For any non-empty $\al\in \hat\cI$, we have
    \[
    |d_{\al}|\leq q^{-|V(\al)|+1} \left(\frac{|a-b|}{n}\right)^{|\al|} f(\al).
    \]
\end{lemma}

\begin{proof}
We proceed by strong induction on $|\al|$. If $|\al|=1$ and $\al\in \hat\cI$, then $\al=\bone_{(1,2)}$. In this case, $d_{\al}=c_{\al}$ by definition. Thus, invoking the expression for $c_{\al}$ in Eq.~\eqref{eq:compute:c}, the conclusion holds for the base case $|\al|=1$.

For the induction step, fix a $\al\in \hat\cI$ and assume the conclusion holds for all $\beta\in \hat\cI$ such that $\be \lneq \al$. Recalling the definition of $d_{\al}$ in \eqref{eq:def:d}, we bound
\begin{equation}\label{eq:bound:d}
|d_{\al}|\leq |c_{\al}|+\sum_{\be\ga\in \hat\cJ \, : \, \be\lneq\al}|u_{\be\ga}|\cdot|M_{\be\ga,\al}|.
\end{equation}
We focus on the sum in \eqref{eq:bound:d}. Since $u_{\emptyset}=0$, we can restrict the sum to $\be\neq \emptyset$. Using Lemma~\ref{lem:bound:M},
\[
\sum_{\be\ga\in \hat\cJ \, : \, \be\lneq\al}|u_{\be\ga}|\cdot|M_{\be\ga,\al}|\leq \sum_{\substack{\be\ga\in \hat\cJ\\\emptyset \lneq \be\lneq\al}}|u_{\be\ga}| \cdot M_{\be\ga,\be}\left(\frac{|a-b|}{n}\right)^{|\al|-|\be|} q^{-|V(\al)|+|V(\be)|}.
\]
For a fixed $\beta\in \hat\cI$ such that $\emptyset\lneq \beta\lneq \al$, our choice of $u_{\be\ga}$ in Eq.~\eqref{eq:sbm:u} implies
\[
\sum_{\ga\in [q]^{V(\be)}}|u_{\be\ga}| \cdot M_{\be\ga,\be}=|d_{\be}|\leq q^{-|V(\be)|+1} \left(\frac{|a-b|}{n}\right)^{|\be|} f(\be)
\]
where the inequality holds by the induction hypothesis. Thus, it follows that
\[
\sum_{\be\ga\in \hat\cJ \, : \, \be\lneq\al}|u_{\be\ga}|\cdot|M_{\be\ga,\al}|\leq q^{-|V(\al)|+1}\left(\frac{|a-b|}{n}\right)^{|\al|} \sum_{\be\in \hat\cI \, : \, \emptyset \lneq \beta\lneq \al} f(\be).
\]
Plugging this estimate into \eqref{eq:bound:d} and using the expression \eqref{eq:compute:c} for $|c_{\al}|$, it follows that
\[
|d_{\al}|
\leq q^{-|V(\al)|+1}\left(\frac{|a-b|}{n}\right)^{|\al|} \left(1+\sum_{\be\in \hat\cI \,:\, \emptyset\lneq \beta\lneq \al} f(\be)\right)=q^{-|V(\al)|+1}\left(\frac{|a-b|}{n}\right)^{|\al|} f(\al),
\]
which concludes the proof.
\end{proof}

By invoking the formula for $M_{\be\ga,\al}$, we can express
\begin{equation}\label{eq:express:sos:M}
\sum_{\ga\in [q]^{V(\al)}}M_{\al\ga,\al}^2=\EE\Bigg[\left(\frac{a}{n}\left(1-\frac{a}{n}\right)\right)^{\ell(\al,\sigma^\star)}\left(\frac{b}{n}\left(1-\frac{b}{n}\right)\right)^{|\al|-\ell(\al,\sigma^\star)}\Bigg],
\end{equation}
where the expectation is over the community labels $\sigma^\star=(\sigma^\star_v)\stackrel{i.i.d.}{\sim}\textnormal{Unif}([q])$. As a result, we establish the following lower bound. Recall the parameter $\xi\equiv \min\{a(1-\frac{a}{n}),b(1-\frac{b}{n})\}/d$.
\begin{lemma}
    For any non-empty $\al\in \hat\cI$,
    \[
    \sum_{\ga\in [q]^{V(\al)}}M_{\al\ga,\al}^2\geq \max\left\{\left(\frac{\dc}{n}\right)^{|\al|}\cdot \xi^{|\al|-|V(\al)|+1}\,,\,\left(\frac{a}{n}\left(1-\frac{a}{n}\right)\right)^{|\al|}\cdot q^{-|V(\al)|+1}\right\}.
    \]
\end{lemma}
\begin{proof}
Recall the expression \eqref{eq:express:sos:M} for $\sum_{\ga\in [q]^{V(\al)}}M_{\al\ga,\al}^2$. On the event that $\sigma^\star_u=\sigma^\star_v$ for all $u,v\in V(\al)$, which happens with probability $q^{-|V(\al)|+1}$, we have $\ell(\al,\sigma^\star)=|\al|$. Thus,
\[
  \sum_{\ga\in [q]^{V(\al)}}M_{\al\ga,\al}^2\geq q^{-|V(\al)|+1}\cdot \left(\frac{a}{n}\left(1-\frac{a}{n}\right)\right)^{|\al|}.
\]
Turning to the other bound, recall that any non-empty $\al \in \hat\cI$ must be connected. Let $\al^\diamond$ be a spanning tree of $\al$. Since $\dc\xi = \min\{a(1-\frac{a}{n}),b(1-\frac{b}{n})\}$, we can lower bound the RHS of~\eqref{eq:express:sos:M} by
\[
\sum_{\ga\in [q]^{V(\al)}}M_{\al\ga,\al}^2\geq (\dc \xi)^{|\al|-|\al^{\diamond}|} \; \EE\Bigg[\left(\frac{a}{n}\left(1-\frac{a}{n}\right)\right)^{\ell(\al^{\diamond},\sigma^\star)}\left(\frac{b}{n}\left(1-\frac{b}{n}\right)\right)^{|\al^{\diamond}|-\ell(\al^{\diamond},\sigma^\star)}\Bigg].
\]
Since $\al^{\diamond}$ is a spanning tree, we have $|\al|-|\al^{\diamond}|=|\al|-|V(\al)|+1$. Moreover, the expectation above can be computed explicitly as follows. Since $\al^{\diamond}$ is a tree and $\sigma^\star$ is a uniform coloring of the vertices, the number of monochromatic edges is distributed according to $\ell(\al^{\diamond},\sigma^\star)\sim\textnormal{Binomial}(|\al^{\diamond}|,1/q)$. Thus, 
\[
\EE\Bigg[\left(\frac{a}{n}\left(1-\frac{a}{n}\right)\right)^{\ell(\al^{\diamond},\sigma^\star)}\left(\frac{b}{n}\left(1-\frac{b}{n}\right)\right)^{|\al^{\diamond}|-\ell(\al^{\diamond},\sigma^\star)}\Bigg]=\left(\frac{1}{q}\cdot \frac{a}{n}\left(1-\frac{a}{n}\right)+\frac{q-1}{q}\cdot \frac{b}{n}\left(1-\frac{b}{n}\right)\right)^{|\al^{\diamond}|},
\]
which equals $(\dc/n)^{|\al^{\diamond}|}$. Therefore, combining the two displays above completes the proof.
\end{proof}

\subsection{Proof of Theorem~\ref{thm:sbm:low:degree}}

We begin with the proof of the first item. Invoking the bound~\eqref{eq:sbm:corr:bound}, 
\[
\begin{split}
\Corr^2_{\leq D}
&\leq \frac{q}{1-1/q} \sum_{\al\in \hat\cI \, : \, 1\leq |\al|\leq D} \frac{d_\al^2}{\sum_{\ga\in [q]^{V(\al)}}M_{\al\ga,\al}^2}\\
&\leq \frac{q}{1-1/q} \sum_{\al\in \hat\cI \, : \, 1\leq |\al|\leq D}\frac{q^{-2|V(\al)|+2} \left(\frac{|a-b|}{n}\right)^{2|\al|} f(\al)^2}{\left(\frac{d_{\circ}}{n}\right)^{|\al|}\xi^{|\al|-|V(\al)|+1}},
\end{split}
\]
where we used Lemma~\ref{lem:d:alpha} and Lemma~\ref{lem:bound:M} in the final inequality.

We will break up the sum according to the parameters $t := |\alpha|$ and $\Delta := t-v+1 \ge 0$ where $v := |V(\alpha)|$. From Lemma~\ref{lem:f:bound}, we have $f(\alpha) \le (2t)^{t-v+1} = (2t)^\Delta$. From Lemma~5.5 in~\cite{SW:24} (and recalling our Remark~\ref{rmk:comparison:sw} on how our conventions compare to that paper), the number of informative $\alpha$ with $|\alpha| = t$ and $|V(\alpha)| = v$ is at most $n^{v-2} (2t)^{5\Delta}$. Also recall that $d_{\circ}\la_{\circ} = \frac{a-b}{q}$. Putting it together, we have
\[
\begin{split}
\Corr^2_{\leq D}
&\leq \frac{q}{1-1/q} \sum_{t=1}^{D}\sum_{\Delta=0}^{t-1} n^{v-2}(2t)^{5\Delta} \cdot \frac{q^{-2v+2} (d_{\circ}^2 \la_{\circ}^2 q^2)^t n^{-2t} (2t)^{2\Delta}}{\left(\frac{d_{\circ}}{n}\right)^t \xi^\Delta}  \\
&= \frac{q}{1-1/q} \sum_{t=1}^{D}\sum_{\Delta=0}^{t-1}\frac{q^{2\Delta} (d_{\circ}^2 \la_{\circ}^2)^{t}(2t)^{7\Delta}}{n^{\Delta+1}d_{\circ}^{t}\xi^{\Delta}}\\
&\leq \frac{q}{1-1/q} \cdot \frac{1}{n}\sum_{t=1}^{D}(d_{\circ}\la_{\circ}^2)^t \sum_{\Delta \geq 0}\left(\frac{q^2(2D)^7}{n\xi}\right)^{\Delta},
\end{split}
\]
which concludes the proof of the first item.

For the second item, we proceed in the same manner except that we use the alternate lower bound of $\sum_{\ga\in [q]^{V(\al)}}M_{\al\ga,\al}^2$ stated in Lemma~\ref{lem:bound:M}. This leads to
\[
\begin{split}
\Corr^2_{\leq D}
&\leq \frac{q}{1-1/q} \sum_{\al\in \hat\cI \, : \, 1\leq |\al|\leq D}\frac{q^{-2|V(\al)|+2} \left(\frac{|a-b|}{n}\right)^{2|\al|} f(\al)^2}{\left(\frac{a}{n}\left(1-\frac{a}{n}\right)\right)^{|\al|}q^{-|V(\al)|+1}}\\
&\leq \frac{q}{1-1/q} \sum_{t=1}^{D}\sum_{\Delta=0}^{t-1}\frac{q^{2\Delta} (d_{\circ}^2 \la_{\circ}^2)^{t}(2t)^{7\Delta}}{n^{\Delta+1}\left(\frac{a(1-\frac{a}{n})}{q}\right)^{t}q^{\Delta}}\\
&\leq \frac{q}{1-1/q} \cdot \frac{1}{n}\sum_{t=1}^{D}(d_{\circ}\la_{\circ}^2)^t\cdot \left(1+\frac{(q-1)b(1-\frac{b}{n})}{a(1-\frac{a}{n})}\right)^t \sum_{\Delta \geq 0}\left(\frac{q(2D)^7}{n}\right)^{\Delta},
\end{split}
\]
where in the last inequality, we used $\frac{a(1-\frac{a}{n})}{q}=d_{\circ}\left(1+\frac{(q-1)b(1-\frac{b}{n})}{a(1-\frac{a}{n})}\right)^{-1}$. This concludes the proof of Theorem~\ref{thm:sbm:low:degree}.

\subsection{Graphon estimation implications -- Corollary \ref{cor: graphon estimation}}\label{subsec: graphon estimation}
\begin{proof}[Proof of Corollary \ref{cor: graphon estimation}]
    Consider $\P_{\xi}$ to be uniform on $[k]$. For $p=(p_1,p_2)$, consider  $M(p) \in \cM_q$, where $M_{i,j}(p)=p_2\one\{z_i=z_j\}+p_1\one\{z_i\neq z_j\}$ if $i\neq j$. Then
    \[
    \begin{split}
    \sup_{M\in \cM_q}\inf_{\wM\in \R_{\leq D}^{n\times n}[Y]} \EE\Big[\ell(\wM,M)\Big]
    &\geq \inf_{\wM\in \R_{\leq D}^{n\times n}[Y]}\EE\Big[\ell\big(\wM,M(p)\big)\Big]\\
    &=\inf_{f\in \R[Y]_{\leq D}}\EE\left[\left(f(Y)-M_{1,2}(p)\right)^2\right].
    \end{split}
    \]
    By a standard fact (see e.g.\ Fact 1.1 in \cite{SW:22}), the RHS equals
    \[
    \begin{split}
    \inf_{f\in \R[Y]_{\leq D}}\EE\left[\left(f(Y)-M_{1,2}(p)\right)^2\right]
    &= \Var{M_{1,2}(p)}\left(1-\sup_{f\in \R[Y]_{\leq D}}\left\{\frac{\Cov{f(Y)} {M_{1,2}(p)}}{\sqrt{\EE[f(Y)^2]}\cdot \sqrt{\EE[M_{1,2}(p)^2]}}\right\}\right)\\
    &=\frac{(p_2-p_1)^2}{q}\left(1-\frac{1}{q}\right)\cdot\left(1-\Corr^2_{\leq D}\right),
    \end{split}
    \]
    where the last equality holds because $M_{1,2}(p)=p_1+(p_2-p_1)\one\{z_i=z_j\}$. By Theorem~\ref{thm: main hard}, as long as $ p_1\leq p_2\leq q p_1$ and 
    \[
    \frac{n(p_2-p_1)^2}{q\big(p_1(1-p_1)+(q-1)p_2(1-p_2)\big)}\leq 1,\quad\quad D\leq n^{c}
    \]
    then $\Corr_{\leq D} \leq C\sqrt{\frac{Dq}{n}}$ where $c, C>0$ depend on $\eps$. Note that we can choose $p_1,p_2$. The choice $p_2=\frac{1}{2}$ and $p_1=\frac{1}{2}-\frac{1}{4}\sqrt{\frac{q^2}{n}}$ satisfies the condition above as long as $2\leq q\leq n^{1/2-\eps}$ and $D\leq n^c$. Therefore,
    \[
      \inf_{\wM\in \R_{\leq D}^{n\times n}[Y]}\sup_{M\in \cM_q} \EE\Big[\ell(\wM,M)\Big]\geq \frac{q}{8n},
    \]
    which concludes the proof.
\end{proof}

\section{Information-theoretic threshold}\label{sec: IT}

\begin{algorithm}[H]\label{alg: inefficient}
    \SetAlgoLined
    Partition the edges $E(G)$ uniformly at random into two sets $E_1$ and $E_2$. \\
    Let $\tau$ be a good partition of $G = (V, E_1)$ as found in \cite[Theorem 9]{BMNN:16}. \\
    \For{$v \in V$}{
    Run depth-1 belief propagation using edge set $E_2$ and initialization $\tau$ to obtain $\widehat\sigma(v)$.
    }
    \caption{Inefficient algorithm}
\end{algorithm}

In this section, we complement our story in the algorithmic regime with an information-theoretic analysis. Following our investigation into the performance of efficient algorithms, it is natural to consider the performance of any possible algorithm on the weak recovery task. We prove Theorem~\ref{thm: IT} and show that the natural scaling for the regime in which exponential-time algorithms succeed is $d\lambda = \Theta(1)$. The results of this section are inspired by the argument of Banks, Moore, Neeman, and Netrapalli~\cite{BMNN:16} in the fixed community setting. Their main idea is as follows. We can now afford a brute-force search over all possible partitions of the vertex set, so it suffices to find a condition under which the partition correlates with the true partition with high probability. It turns out that checking the number of edges between different parts of the partition is a sufficient condition for this task. Analyzing these edge counts yields the following result in~\cite{BMNN:16}.
\begin{theorem}[\cite{BMNN:16}, Theorem~9]\label{thm: weak IT}
    When
    \[ d > \frac{2q \log q}{(1+(q-1)\lambda)\log(1+(q-1)\lambda) + (q-1)(1-\lambda)\log(1-\lambda)} \]
    we can recover $\beta$ fraction of the vertices (inefficiently) where $\beta$ is the smallest solution to 
    \[ d = \frac{2(h(\beta) + (1-\beta)\log q)}{-\lambda\log \beta + (1-\lambda)\log\frac{1-\lambda}{1-\beta\lambda}}. \]
\end{theorem}
We unpack the implications of this bound in our regime, asymptotically as $n\to\infty$. It is equivalent to recovering a $\beta$ fraction of labels when $d \gtrsim \frac{2}{\lambda}$, where $\beta$ satisfies the following:
\[ d = \frac{2(1-\beta)}{-\lambda\frac{\log \beta}{\log q}} \implies \beta = (1+o(1))\left(\frac{1}{q}\right)^{\frac{2}{d\lambda}}. \]
The idea is to boost this via an application of belief propagation to get $\delta$-recovery. 

\begin{proof}[Proof of Theorem~\ref{thm: IT}]
We begin with the achievability result. We first analyze the effect of running one step of belief propagation. Recall that the belief propagation recursion on a tree is given by the following (see e.g.\ \cite{S:11}):
\[ X_\rho(j) = \frac{\prod_{i=1}^d [1+\lambda q(X_{u_i}(j) - \frac{1}{q})]}{\sum_{k=1}^q \prod_{i=1}^d [1+\lambda q(X_{u_i}(k) - \frac{1}{q})]}. \]
Here $X_v(j)$ is the probability that vertex $v$ is in community $j$. The Banks--Moore--Neeman--Netrapalli algorithm outputs an assignment of communities satisfying $X_v(j) = \beta$ for $j = \sigma^\star_v$ and $j = \frac{1-\beta}{q-1}$ otherwise. For each vertex, we expect $\frac{d(1 +(q-1)\lambda)}{q}$ neighbors of the matching community, and $\frac{d(1-\lambda)}{q}$ neighbors of each other community. Recall that we are in a regime where $d\lambda \geq 14$ and $\frac{d}{q} = O(\frac{1}{q})$. Thus with positive probability we will have at least $\frac{d\lambda}{2}$ neighbors of the same community, and at most 2 neighbors of any other community. On this event, the products in the recursion satisfy 
\[ \prod_{i=1}^d [1+\lambda q(X_{u_i}(j) - \frac{1}{q})] \geq (1-\lambda \beta)^d\left(\frac{1+\lambda q\beta}{1-\lambda \beta}\right)^{\frac{d\lambda}{2}} \]
when $j = \sigma^\star_\rho$ and
\[ \prod_{i=1}^d [1+\lambda q(X_{u_i}(j) - \frac{1}{q})] \leq (1-\lambda \beta)^d\left(\frac{1+\lambda q\beta}{1-\lambda \beta}\right)^2 \]
when $j \neq \sigma^\star_\rho$. Notice that for our value of $\beta$, we have
\[ \frac{1+\lambda q\beta}{1-\lambda \beta} \approx \lambda q \beta = \lambda q^{1-\frac{2}{d\lambda}}. \]
Thus, we have that 
\[ X_\rho(\sigma^\star_\rho) \geq \frac{\left(\lambda q^{1-\frac{2}{d\lambda}}\right)^{\frac{d\lambda}{2}}}{\left(\lambda q^{1-\frac{2}{d\lambda}}\right)^{\frac{d\lambda}{2}} + q\left(\lambda q^{1-\frac{2}{d\lambda}}\right)^2} \asymp \frac{q^{\frac{d\lambda}{2}-1}}{q^{\frac{d\lambda}{2}-1} + q^{3-\frac{4}{d\lambda}}}. \]
When $d\lambda > 7$ we have $\frac{d\lambda}{2}-1 > 3-\frac{4}{d\lambda}$ so in particular $X_\rho(\sigma^\star_\rho) = 1 - o(1)$. Thus, conditioned on this good event in the neighborhood of a vertex, one iteration of belief propagation recovers the root with high probability. We now use this to obtain an $\epsilon$ alignment partition of the vertices.

Suppose $d\lambda > 14$. Partition the edges of the graph uniformly at random into two equal sized sets $E_1$ and $E_2$. Then conditioned on the true community labeling, each set is an $\SBM(n,q,d,\lambda)$ with $d\lambda > 7$. Applying Theorem~\ref{thm: weak IT} on the graph with edges $E_1$, we perform a brute-force search to find a good partition with alignment $\beta \geq \left(\frac{1}{q}\right)^{\frac{2}{d\lambda/2}}$. Now for each vertex $v$, the above calculation shows that the BP estimator is a function of its $E_2$-neighborhood that gives the correct label with probability at least $\epsilon > 0$. By Lemma~\ref{lemma: coupling}, we can partition the vertices into groups of size $\log n$ such that the estimators within each group are approximately independent. This implies that the total fraction of labels recovered correctly concentrates around $\epsilon$. 

Part 2 of the theorem --- the impossibility result --- follows from a percolation style argument. We first show that we may reduce to the case where only the within-community edges are revealed to us. Consider generating a sample from $G \sim \SBM(n,q,d,\lambda)$ as follows. Assign each vertex $v$ an i.i.d.\ uniform $\sigma^\star_v \in [q]$ community label as before. For each pair of vertices $u,v$ such that $\sigma^\star_u = \sigma^\star_v$ connect them with an edge independently with probability $\frac{qd\lambda}{n} \cdot \frac{1}{1-\frac{d(1-\lambda)}{n}} = (1+o(1))\frac{qd\lambda}{n}$ to obtain a graph $\widetilde{G}$. Finally, let $G = \widetilde{G} \cup \widehat{G}$ where $\widehat{G}$ is the Erd\H{o}s--R\'enyi graph $G\left(n, \frac{d(1-\lambda)}{n}\right)$ independent of $\widetilde{G}$. We claim that $\delta$-recovery is possible in $\SBM(n,q,d,\lambda)$ only if recovery is possible using $\widetilde{G}$. Indeed, given a sample of $\widetilde{G}$ we can simply sample $\widehat{G}$ and perform recovery from the block model. Thus, it suffices to show that $\delta$-recovery is impossible from $\widetilde{G}$. 

$\widetilde{G}$ is the union of $q$ copies of $G\left(\frac{n}{q}, (1+o(1))\frac{qd\lambda}{n}\right)$, each of which is sub-critical since $d\lambda < 1$.  Classic theory of random graphs, present even in the initial work~\cite{ER:60}, implies that each component of $\widetilde{G}$ has size $O(\log \frac{n}{q})$ with high probability. Now consider any algorithm which takes as input $\widetilde{G}$ and outputs an assignment $\widehat{\sigma} \in [q]^V$. Let $\widehat{V}_i = \{v: \widehat{\sigma}_v = i\}$ be the $i$th community of output and $V^\star_i = \{v: \sigma^\star_v = i\}$ be the true $i$th community. We may assume that each $\widehat{V}_i$ consists of a union of components, since the information of the input graph exactly informs us that connected components are in a single community. Suppose we are even informed of one component from each community --- this corresponds to selecting a permutation from the maximum in Definition~\ref{def: Align}. By symmetry, the remaining individual components are exchangeable with respect to the communities. Thus, each remaining component is assigned correctly with probability $\frac{1}{q}$. The expected alignment is still at most 
\[ \frac{1}{n}\left(O(q\log \frac{n}{q}) + \frac{1}{q}(n-O(q\log \frac{n}{q}))\right) = O\left(\frac{q\log \frac{n}{q}}{n} + \frac{1}{q}\right) = o(1) \]
so $\delta$-recovery is impossible from $\widetilde{G}$ for any $\delta > 0$ and $q \ll n$. 
\end{proof}

\section{Detection via triangle counts}\label{sec: detection}
To conclude, we consider the problem of detection --- that is distinguishing a graph drawn from the SBM from one drawn from an Erd\H{o}s--R\'enyi distribution with the same average degree. Our main result is that the detection problem is efficiently solvable whenever $q$ is increasing and $\lambda \neq 0$. 

While the conditions of the proposition may be difficult to digest at first, consider the setting in which $d, \lambda$ are fixed while $n \to \infty$. Our result says that in this regime, as long as $q \to \infty$ and $\lambda \neq 0$ the triangle count solves the detection problem, yielding Theorem~\ref{thm: detection} as a special case. 
\begin{proposition}
    Suppose $q\lambda \to \infty$, $d^3q\lambda^3 \to \infty$ and $nd\lambda \to \infty$. Then as long as $q\lambda^3$ is bounded away from 1 there exists a polynomial-time algorithm that distinguishes between $G(n, \frac{d}{n})$ and $\SBM(n, q, d, \lambda)$ with probability $1-o(1)$. 
\end{proposition}
\begin{proof}
The algorithm is based on counting the number of triangles in the sampled graph. Recall that in an Erd\H{o}s--R\'enyi random graph $G(n, \frac{d}{n})$ the number of triangles is distributed approximately as $\mathrm{Pois}(\frac{d^3}{6})$ when $d = \Theta(1)$ and satisfies a central limit theorem when $d = \omega(1)$~\cite{R:88}. We analyze the distribution of triangles under $\mathrm{SBM}(n,q,d,\lambda)$. 
\begin{align*}
    \E{N_{K_3}(\mathrm{SBM}(n,q,d,\lambda))} &= \binom{n}{3}\left(\frac{1}{q^2}\cdot \frac{a^3}{n^3} + \frac{(q-1)(q-2)}{q^2}\cdot \frac{b^3}{n^3} + \frac{3(q-1)}{q^2}\cdot\frac{ab^2}{n^3}\right) \\
    &= (1+o(1))(\frac{a^3}{6q^2} + \frac{ab^2}{2q} + \frac{b^3}{6}) \\
    &= (1+o(1)) d^3\left( \frac{(1+q\lambda)^3}{6q^2} + \frac{1+q\lambda}{q} + \frac{1}{6} \right).
\end{align*}
If $q\lambda \to \infty$ then asymptotically we have $\frac{d^3}{6}q\lambda^3$. In this regime, we can compute the variance as 
\begin{align*}
    \Var{N_{K_3}(\mathrm{SBM}(n,q,d,\lambda))} &= \E{N_{K_3}(\mathrm{SBM}(n,q,d,\lambda))} + (1+o(1))4\binom{n}{4}\cdot \frac{1}{q^3}\left(\left(\frac{a}{n}\right)^5 - \left(\frac{a}{n}\right)^6\right) \\
    &= (1+o(1))\left(\frac{d^3}{6}q\lambda^3 + \frac{(dq\lambda)^5}{6nq^3} \right) \\
    &= (1+o(1))\left(\frac{d^3}{6}q\lambda^3 + \frac{d^5q^2\lambda^5}{6n} \right).
\end{align*}
By Chebyshev's inequality, 
\begin{align*}
    \Prob{\abs{N_{K_3}(\mathrm{SBM}(n,q,d,\lambda)) - \frac{d^3}{6}q\lambda^3} > \epsilon d^3 q\lambda^3} &\leq (1+o(1))\left(\frac{d^3}{6}q\lambda^3 + \frac{d^5q^2\lambda^5}{6n} \right) \cdot \frac{1}{\epsilon^2d^6q^2\lambda^6} \\
    &= (1+o(1))\left( \frac{1}{6\epsilon^2d^3q\lambda^3} + \frac{1}{6\epsilon^2nd\lambda} \right).
\end{align*}
Thus, as long as $d^3q\lambda^3 \to \infty$ and $nd\lambda \to \infty$ we have concentration of the triangle count around its expectation. The condition $q\lambda^3$ being bounded away from 1 ensures that the two models have different expectations. This implies that counting triangles distinguishes the two distributions with high probability. 
\end{proof}

\bibliographystyle{amsplain0}
\bibliography{SBM}

\end{document}